\documentclass[review]{elsarticle}

\usepackage{hyperref}

\usepackage[margin=2.5cm]{geometry}
\usepackage{url}        	
\usepackage{booktabs}   	
\usepackage{amsfonts}   	
\usepackage{nicefrac}   	
\usepackage{microtype}  	
\usepackage{lipsum}
\usepackage{fancyhdr}   	
\usepackage{bm}
\usepackage{bbold}
\usepackage{graphicx}   	
\usepackage{subcaption}
\usepackage[inline]{enumitem}
\usepackage{diagbox}
\usepackage{setspace} 
\usepackage{cancel}
\usepackage{mathtools}
\usepackage{comment}
\usepackage{amssymb}
\usepackage{amsthm}
\usepackage{amsfonts}
\usepackage{bm}
\usepackage{amsbsy}
\usepackage{amsmath}
\usepackage{float}
\usepackage{mathrsfs}
\usepackage[nameinlink]{cleveref}
\usepackage{color}
\usepackage{verbatim}
\usepackage[linesnumbered, ruled,vlined]{algorithm2e}
\usepackage{multirow}
\usepackage{booktabs}
\usepackage{array}
\usepackage{adjustbox}
\usepackage{longtable}
\usepackage{pdflscape}
\usepackage[numbers]{natbib}

\SetKwComment{Comment}{$\triangleright$\ }{}
\SetCommentSty{}
\SetKw{Continue}{continue}
\SetKw{Break}{break}

\newcommand{\bR}{\mathbb{R}}
\newcommand{\cA}{{\mathcal A}}

\newcommand{\cF}{{\mathcal F}}
\newcommand{\cE}{{\mathcal E}}
\newcommand{\cEI}{{\mathcal {EI}}}

\newcommand{\cP}{{\mathcal P}}

\newcommand{\cV}{{\mathcal V}}

\newcommand{\st}{{\text{ s.t. }}}

\newcommand{\dlt}{{\delta}}
\newcommand{\problem}{CSCP$_\dlt$}
\newcommand{\dfproblem}{discrete facility SCP$_\dlt$}

\newcommand{\ncalgo}{\text{nodeCover}}
\newcommand{\malgo}{\text{mutual}}

\newcommand{\todo}[1]{{\color{red}#1}}
\newcommand{\bO}[1]{{\mathcal{O}(#1)}}
\newcommand{\dash}{:}
\newcommand{\tabledefline}[2]{\multicolumn{1}{l}{\rlap{#1\ \dash\ #2}}\\}

\DeclareMathOperator{\prev}{prev}

\DeclarePairedDelimiter\ceil{\lceil}{\rceil}

\DeclareMathOperator*{\argmin}{arg\,min}

\theoremstyle{definition}
\newtheorem{definition}{Definition}[section]
\newtheorem{observation}{Observation}[section]
\newtheorem{theorem}{Theorem}[section]
\newtheorem{lemma}{Lemma}[section]
\newtheorem{proposition}{Proposition}[section]

\newtheorem{assumption}{Assumption}[section]
\newtheorem{remark}{Remark}[section]

\begin{document}

\begin{frontmatter}

\title{Continuous covering on networks: Improved mixed integer programming formulations}

\author[1,2,address]{Mercedes Pelegrín}
\ead{pelegringarcia@lix.polytechnique.fr}

\author[1,address]{Liding Xu\corref{mycorrespondingauthor}}
\cortext[mycorrespondingauthor]{Corresponding author}
\ead{liding.xu@polytechnique.edu}

\address[address]{Laboratoire d'Informatique (LIX), École Polytechnique, 91128 Palaiseau, France}

\fntext[fn1]{The authors have contributed equally.}
\fntext[fn2]{Mercedes Pelegrín was with Laboratorire d'Informatique de l'École Polytechnique. She is currently with EMEA Analytics at FICO.}
\begin{abstract}
Covering problems are well-studied in the domain of Operations Research, and, more specifically, in Location Science.
When the location space is a network, the most frequent assumption is to consider the candidate facility locations, the points to be covered, or both, to be finite sets. In this work, we study the set-covering location problem when both candidate locations and demand points are continuous on a network. This variant has received little attention, and the scarce existing approaches have focused on particular cases, such as tree networks and integer covering radius. Here we study the general problem and present a Mixed Integer Linear Programming formulation (MILP) for networks with edge lengths no greater than the covering radius. The model does not lose generality, as any edge not satisfying this condition can be partitioned into subedges of appropriate lengths without changing the problem. We propose a preprocessing  algorithm to reduce the size of the MILP, and devise tight big-M constants and valid inequalities to strengthen our formulations. Moreover, a second MILP is proposed, which admits edge lengths greater than the covering radius. As opposed to existing formulations of the problem (including the first MILP proposed herein), the number of variables and constraints of this second model does not depend on the lengths of the network's edges. This second model represents a scalable approach that particularly suits real-world networks, whose edges are usually greater than the covering radius. Our computational experiments show the strengths and limitations of our exact approach to both real-world and random networks. Our formulations are also tested against an existing exact method.
\end{abstract}

\begin{keyword}
Continuous Facility Location \sep Location on Networks \sep Set-Covering Location Problem \sep Mixed Integer Programming
\end{keyword}

\end{frontmatter}

\section{Introduction}
Covering in Operations Research refers to the optimization problem of deciding the location of facilities to ``cover" the points of the so-called demand set, which should fall within the radius coverage of at least one of the installed facilities.
This classic problem finds applications in many different domains, including health care \cite{health}, surveillance of transport networks \cite{Gusev20}, computer networks security \cite{worm}, crane location for construction \cite{crane}, military evacuation systems \cite{military}, homeland defense \cite{homeland}, and urban air mobility \cite{liding}.

Covering problems have taken many forms in the literature. A rough classification distinguishes between \emph{maximal covering} location and \emph{set-covering} location problems. The former aims at maximizing the covered demand with a fixed number of facilities (see e.g. \cite{Church74}), while the latter seeks to minimize the number of installed facilities to cover all the demand (see e.g. \cite{Toregas72}). In these classic works \cite{Church74,Toregas72}, the problem is defined on a network and both demand points and candidate facility locations are at nodes. Most of the variants of network covering studied afterward consider at least one of these two sets to be finite, see the reviews \cite{chapter-plastria,chapter-marin} and the references therein. However, this assumption corresponds to ideal but usually unrealistic scenarios (the reader is referred to the real applications of the above paragraph). As an example, in the eVTOLs safety landing site location problem \cite{liding}, a set of emergency landing sites has to be installed on the traffic network in such a way that any point of the same is covered. Something similar happens for the location of ambulance bases in rural areas studied in \cite{health}.
Some works addressing network covering with continuous sets of both candidate locations and demand points are \cite{Carrizosa16,Kalsics,Baldomero}, for maximal covering, and \cite{Gurevich84,Hamacher20,Hartmann21}, for set-covering. We focus on the latter variant, which we call the continuous set-covering problem.

Gurevich et al. \cite{Gurevich84} presented an algorithm to compute an optimal continuous set-covering when the covering radius and the edge's lengths are natural numbers. This algorithm is polynomial time for the class of networks satisfying that every non-separable component is either an edge, a simple cycle, or a simple cycle with one chord, that is, for ``almost tree'' networks. More recently, Fr\"ohlich et al.  \cite{Hamacher20} also studied the same version of the continuous set-covering with natural numbers. The authors presented three different approaches to solve the problem, including a Mixed Integer Linear Programming (MILP) formulation. On the other hand, Hartmann et al. \cite{Hartmann21} focused on the computational complexity of the continuous set-covering for general covering radii. They proved that, when all edges have unit length, the continuous set-covering is polynomially solvable if the covering radius is a unit fraction, and is NP-hard otherwise.

We can now formally state our problem. Consider an undirected connected network $N=(V,E,l)$, where $l:E \to \bR_+$ is the edges' length function. We will denote $l_e:=l(e)$ the length of $e$. The continuum of points on all edges and nodes of $N$ is denoted with $C(N)$. The distance function $d(\cdot,\cdot)$ defines the distance between two points, which coincides with the length of the shortest path in $C(N)$ connecting them.  Given $\dlt > 0$, a point $p\in C(N)$ is said to $\dlt$-cover $p'\in C(N)$  (respectively, $p'$ $\dlt$-covers $p$) if $d(p,p')\le \dlt$ holds. The parameter $\dlt$ is called the covering radius. The continuous $\dlt$-covering location problem on $N$ is to find a set of facility locations in $C(N)$ of minimum cardinality that $\dlt$-covers the whole network, and is formally stated next.
\begin{definition}[Continuous Set-Covering Problem (\problem)]\label{def.problem}
The Continuous Set-Covering Problem on a network $N$ can be expressed as the following optimization problem:
\begin{equation}
\label{cfl}
	\min\Big\{|\cP|: \: \cP=\{p_i\}_{p_i\in C(N)} \textup{ and } \forall p\in C(N), \: \exists p_i\in \cP \st d(p,p_i)\leq \dlt\Big\}.
\end{equation}
A set $\cP$ satisfying the condition within \eqref{cfl} is called a $\dlt$-cover of $N$, while $\cP^*$ minimizing \eqref{cfl} is a minimum $\dlt$-cover.
\end{definition}
The set $\mathcal{P}$ in Definition \ref{def.problem} can represent the locations of ambulance bases  \cite{health}, surveillance cameras \cite{Gusev20}, routing servers in a network of computers \cite{worm}, cranes for construction \cite{crane}, aerial military medical evacuation facilities \cite{military}, aircraft alert sites for homeland defense \cite{homeland}, or eVTOL safety landing sites in an urban area \cite{liding}. 

The \problem\ is known to be NP-hard, see \cite{Hartmann21}.  Due to the continuous nature of \problem, there is an infinite number of candidate  locations. Previous works reduce \problem\ to a tractable set covering problem by discretization. A first observation is that typical simplifications proposed in other related studies are not valid for the \problem.

\subsection{Discretization methods}
\label{sub.disc}

\begin{figure}
	\centering
	\includegraphics[scale=1.4]{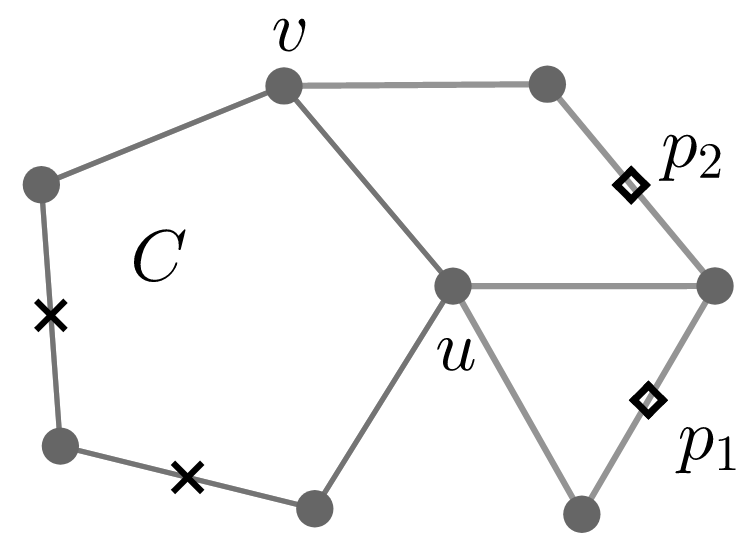}
	\caption{Two cycle coverage points with respect to a cycle $C$ of five nodes}
	\label{fig:cycle-cover}
\end{figure}
Discretization methods identify finite dominating sets (FDS), which are finite subsets of candidate locations guaranteed to contain an optimal solution. 
From the literature, we know at least three FDS for related variants of the \problem, which rely on different assumptions on the network and the covering radius. Here, we review and compare these FDS. Since we aim to propose a general exact algorithm, we show that it may not be viable to extend discretization methods to solve the \problem\ for general networks and real radii.

First, Church and Meadows \cite{Church79} studied the problem with demand at nodes, and identified the following points:
$$NIP:=\{p\in C(N):\: d(p,v)=\dlt \:\textup{ for some }v\in V \}.$$
The authors proved that $FDS_1:=V\cup NIP$ is an FDS for the network set-covering problem when the set of demand points is $V$ and that of candidate locations is $C(N)$.

Secondly, Gurevich et al. \cite{Gurevich84} studied the continuous set-covering problem  when the covering radius and the edge’s lengths are natural numbers. They presented an FDS for the case of all edge lengths being one, which can be easily extended to the case of general edge lengths (see \cite{Hamacher20}),
$$FDS_2:=\big\{p\in C(N):\: d(p,v)=\frac{i}{2\cdot l_e} \:\textup{ for some }e\in E\textup{ and } v\in e;\: i=0,\ldots,2\cdot l_e\big\}.$$
Note that $FDS_2$ depends on the edge's length.

Lastly,  Fr\"ohlich et al. \cite{Hamacher20} proposed a different FDS for the same version of the continuous set-covering with natural numbers. The authors defined the following set of cycle coverage points:
$$CCP:=\left\{p\in C(N):\: d(p,C):=\min_{y\in C}\{d(p,y)\}=\left(\dlt-\frac{l_C}{2}\right)\textup{ mod } \dlt, \:p\notin C,\:\textup{ for a simple cycle }C\subseteq C(N)\right\},$$
where $l_C$ is the total length of the cycle $C$. Suppose that a cycle $C$ is covered by a set of facilities. A cycle coverage point is the furthest point where a facility that contributes to cover $C$ can be moved without compromising the coverage of the cycle (if the rest of the facilities remain unchanged).

Figure \ref{fig:cycle-cover} illustrates this idea. In the depicted example, all edges have unit length and $\dlt=2$. The figure depicts $p_1$ and $p_2$, which are CPP with respect to the cycle $C$ of five nodes. Note that $d(p_1,C)=d(p_2,C)=1.5(=(2-5/2)\textup{ mod } 2)$. Figure \ref{fig:cycle-cover} also depicts two locations in $C$ (marked with symbols 'x'), which correspond to two possible feasible locations for the remaining facility needed to cover $C$ (note that the one at the bottom only yields a covering of the cycle if $p_2$ is located, while the other one together with either $p_1$ or $p_2$ can completely cover $C$).

The authors of \cite{Hamacher20} gave the following recursive definition of an FDS for the problem with natural numbers:
$$S_1:=V\cup NIP \cup CCP;$$
$$S_{j+1}:=S_j\cup\{p\in C(N):\: d(p,y)=\dlt \textup{ for some } y\in bd(\cA(S_j)) \};$$
$$FDS_3:=S_{|J|},$$
where $bd(\cA(S_j))$ is the boundary of the area covered by $S_j$, and  $J\subseteq E$ is the subset of edges to be covered. As the author explained themselves, $FDS_3\subseteq FDS_2$. However, the cardinality of $FDS_3$ may be exponential in the input size, as the number of cycles in a network is in general exponential.

\begin{figure}
	\centering
	\includegraphics{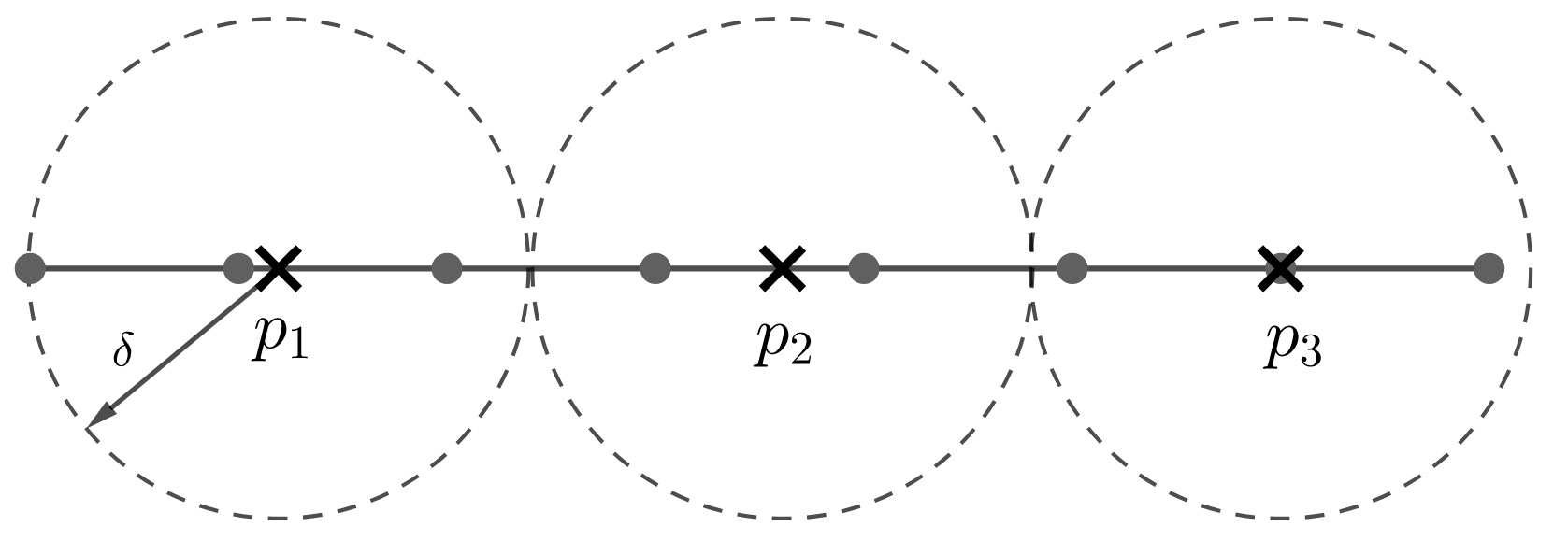}
	\caption{An instance of \problem\ such that not all facilities in $\mathcal{P}^*$ are at a distance $\dlt$ from some node}
	\label{fig:exDiscrete}
\end{figure}
The example depicted in Figure \ref{fig:exDiscrete} illustrates that none of $FDS_1$ and $FDS_2$ are FDS for the \problem. A similar observation was already presented in \cite{Kalsics} for a related problem. The figure shows eight nodes on a path, where all edges have equal lengths. If $l_e=1$ for all $e\in E$ and $\dlt=1.2$, $\cP:=\{p_1,p_2,p_3\}$ is an optimal $\dlt$-cover.

It can be easily observed that there is no optimal solution in which $p_2$ is placed either at a node or at a distance $\dlt$ from some of the eight nodes, which shows that $FDS_1$ is not a valid FDS. On the other hand, it is also easy to check that there is not a feasible solution with the three facilities located either at nodes or middle-points of edges, which proves that $FDS_2$ is also not a valid FDS. As opposed to $FDS_2$, the assumption of the edge lengths and coverage radius being natural numbers is not fundamental in the definition of $FDS_3$. Conversely, $FDS_3$ is based on the idea of identifying those points at the ``boundaries'' of coverage areas, i.e., those delimiting the transition from covering/not covering a specific part of the network.

Consequently, $FDS_3$  could be extended to the general \problem. However, such an extension potentially yields sets with many more candidates, due to the recursive construction of $FDS_3$ based on the distance function. Note that, if $\dlt$ and the edge lengths are natural numbers, $FDS_3$ only contains points of the set
$$INT:=\{p\in C(N): d(p,v)\textup{ is integer or half-integer for some } v\in V\}.$$
Indeed, if $\dlt\in \mathbb{N}$, $NIP\subseteq INT$ is clear; $CCP\subseteq INT$ holds since the operation $(\dlt-{l_C}/{2})\textup{ mod } \dlt$ only yields half integers; and $FDS_3\subseteq INT$ then easily follows by definition. However, if $\dlt\in \mathbb{R}$, the locations of the points in $FDS_3$ are a priori undetermined, and its cardinality increases. Take the same example depicted by Figure \ref{fig:cycle-cover}. If $\dlt=2.1$ (i.e. we increase $\dlt$ just by 0.1), the CCP with respect to the cycle $C$ of the Figure increases from two to four points.

\subsection{Contribution}
As opposed to discretization methods, we directly tackle the \problem\ for general networks and real radii.
Our main contribution is an exact integer programming approach for the \problem, together with tailored algorithms and strategies to tackle it. Even if this problem has been known for decades, surprisingly, only a few partial results are known for some special cases and sub-classes of networks. To the best of our knowledge, only one MILP model \cite{Hamacher20} has been proposed so far which can address the general \problem. Such a model can be applied to any network whose edges do not measure more than the covering radius. This condition does not restrict the applicability of the MILP in \cite{Hamacher20}, as any network can be transformed into an equivalent one that satisfies it.

Here, we present an enhanced MILP formulation that relies on the same assumption as that in \cite{Hamacher20}, but whose numbers of constraints and variables have smaller order of magnitude. In addition, preprocessing strategies to reduce the number of variables of the model are studied, and tailored algorithms are presented. Approaches to strengthen this formulation are also presented, including big-M constants tightening and valid inequalities. The valid inequalities are constraints that reduce the feasible space without removing model solutions.

The introduction of a second MILP, which is scalable concerning the edge's lengths completes the main contributions of this work. This second MILP is an adaptation of the first one we propose, with the difference that it does not require all edge lengths to be smaller than the covering radius. Finally, our computational experiments prove that the MILP model in \cite{Hamacher20}  is not scalable. On the other hand, the preprocessing technique drastically reduces the size of the first model proposed herein.

Finally, we show in the experiments that the second model we propose is superior to both  the model from \cite{Hamacher20} and our first model, in terms of the solution quality and solving time.

In the proposed setting, both the candidate facility locations and the demand points are continuous sets (in particular, they coincide with $C(N)$). The problem could be defined for a subset of demand edges, $J\subseteq E$, and/or a subset of candidate locations $H\subseteq E$. The theoretical results, model, and methods described in this work apply to such cases, after straightforward adaptation.

\subsection{Related works}

Facility location and set covering problems have many variants and applications in operations research and management science. Related literature to this work is vast; here we review a selection of works related to \problem. In \cite{sadigh2010mixed}, the model allows that an edge is covered jointly by two facilities. In \cite{schbel2019}, the authors presented a unified vision of the common characteristics of facility location problems in a continuous space. In \cite{puerto2018extensive}, the authors summarized the research progress in facility location problems on networks. A recent survey \cite{wang2021emergency} provided a comprehensive overview of  emergency facility location problems in logistics, including mathematical models, applications, and the commonly used solution methods.  One of the most distinguishing features of variants of the maximal covering location problem is the solution space: continuous \cite{chapter-plastria, bansal2017planar}, discrete \cite{chapter-marin, cordeau2019benders}, or on networks \cite{bucarey2022benders,berman2016covering}. Especially in \cite{berman2016covering}, the authors  introduced the maximal covering location problem with edge demand. In \cite{baldomero2022upgrading}, the authors studied the upgrading version of the maximal covering location problem with edge length modifications on networks. A related problem that has been recently studied is the obnoxious facility location problem \cite{drezner2018weber}. It aims at locating undesirable facilities that have a negative impact on communities. The most common  objective is to maximize the shortest distance to the closest facility, and the problem has various variants featuring multiple facilities on the plane \cite{drezner2019planar,kalczynski2021obnoxious}, $p$-median objective \cite{kalczynski2022obnoxious}, or edge demand on networks \cite{berman2016covering}. We refer to \cite{church2022review} for a recent review on the obnoxious facility location problem. For more related works, we refer to \cite{akgun2015risk,paul2017multiobjective,fukunaga2016covering,berman2011big,hudec1994confined}.

\subsection{Outline}
The rest of the paper is organized as follows. Section \ref{sec:prelim} presents useful notation and the theoretical development upon which our model is built. Then, our first MILP model is introduced in Section \ref{sec:milp}, while strategies to strengthen this model are described in the next section. The network processing algorithms that complement our MILP are detailed in Section \ref{sec:prepro}. A second MILP model, which we call reduced formulation and is a modification of the first MILP, is presented in Section \ref{sec:reduced-milp}. Finally, Section \ref{sec:compu} describes our computational experiments, and reports and analyzes the obtained results. Section \ref{sec:conclu} closes the paper with some conclusions.

\section{Covering characterization}\label{sec:prelim}
This section presents several notation, definitions, assumptions, observations, and results related to the \problem. On the one hand, Proposition \ref{thm.cover} gives a characterization of the $\dlt$-covers of a network, which is based on the individual coverage of each edge of the network. Then, this result is refined to obtain a second necessary and sufficient covering condition in Proposition \ref{cora.covercond}. It distinguishes between two alternative possibilities for covering each edge, namely complete or partial, and will be useful for our MILP formulation and methods. The rest of the section is oriented to characterize the so-called partial and complete covers. The idea of these sets is to delimit the areas of the network where a facility, if placed, would completely cover a given edge, and those where a facility would reach the edge (but maybe not completely cover it).

We first introduce some related notation, definitions and assumptions.
We assume that $V$ is totally ordered by the binary relation $\preceq$. Every edge $e \in E$ has a unique representation, $e = (v_a, v_b)$, where $v_a, v_b \in V$, and $v_a \preceq v_b$. From now on, we take $e=(v_a, v_b)$ indifferently as a continuum in $C(N)$ or as an edge ending at  $v_a, v_b$. We extend the edges' length function to  $l: C(N) \to \bR_+$ as a length measure on the continuum of points.   For two points $p, p' \in C(N)$, we denote by $\Pi(p,p') \subseteq 2^{C(N)}$ and $\Pi^*(p,p') \subseteq \Pi(p,p')$ the set of paths and shortest paths, respectively, connecting $p$ and $p'$. Any path $\pi\in \Pi(p,p')$ is indifferently treated as a continuum in $C(N)$, then $l_{\pi}:=l(\pi)$ is the length of $\pi$.  The distance between $p$ and  $p'$, $d(p, p')$, is the length of a shortest path connecting them:
$$d(p,p') := \min\{l_{\pi}: \:\pi\in \Pi(p,p')\}=l_{\pi^*} \textup{ for any } \pi^*\in \Pi^*(p,p').$$
In particular, if $p$ and $p'$ belong to the same edge, we denote by $l(p,p')$ the length of the unique path in that edge connecting them.
We work under the following assumption:
\begin{assumption} \label{assumption}
$\dlt \ge l_e$ for all $e \in E$.
\end{assumption}
If Assumption \ref{assumption} did not hold, we could consider a set $I \subset  2^{C(N)}$, that would contain, for each $e = (v_a, v_b) \in E$, the following continuum sets of points (segments):
\begin{itemize}
	\item If $\dlt \ge l_e$, $e \in I$;
	\item If $\dlt < l_e$, let  $n:= \ceil{\frac{l_e}{\dlt}}+1$. We define $v_1:= v_a$, $v_n:= v_b$ and $v_2,\hdots, v_{n-1}\in e$ such that $l(v_a,v_i) = (i -1) \frac{l_e}{n-1}$ for $i=2,\ldots, n-1$. Then, $(v_{i}, v_{i+1}) \in I$ for all $i=1,\ldots,n-1$.
\end{itemize}
We consider $N'=(V',E'=I)$, where $V'$ contains the endpoints of $I$. The new network $N'$ satisfies that $\dlt \ge l_e$ for all $e\in E'$, and it is isomorphic to $N$ with respect to the length function. Indeed, since $N'$ is obtained by subdividing edges in $N$,  $C(N)=C(N')$ and a set of points $\dlt$-covers $N$  if and only it $\dlt$-covers $N'$. Therefore, Assumption \ref{assumption} always holds after the network $N$ is transformed into $N'$ (via a preprocessing step). Such transformation yields a network with more nodes and edges, which has a direct impact on the size of optimization models. In \Cref{sec:reduced-milp}, we present a model that avoids this effect.  

\subsection{Observations}
In the following, we give several observations of optimal $\dlt$-covers, which guide our quantitative analysis of covering conditions and the resulting MILP model of Section \ref{sec:milp}.

\begin{observation}
For an edge $e=(v_a, v_b) \in E$ and a point $p \in C(N)$, one of the following cases holds:
\begin{enumerate}
	\item   $p$ cannot $\delta$-cover any point in $e$;
	\item $p$ can $\delta$-cover the whole $e$;
	\item $p$ can $\delta$-cover a continuous portion of $e$ containing either $v_a$ or $v_b$;
	\item $p$ can $\delta$-cover two continuous portions of $e$, which do not intersect, each contains either $v_a$ or $v_b$.
\end{enumerate}
\end{observation}

\begin{observation}{(a similar statement was proven in \cite{Hamacher20})}
There exists an optimal $\dlt$-cover that satisfies:
\begin{enumerate}
	\item[i)] Each edge $e \in E$ has at most two facilities (due to Assumption \ref{assumption});
	\item[ii)] If there are two facilities in the edge $e$, we can assume without loss of generality that they are located at end nodes $v_a$, $v_b$ (this follows from i) and  Assumption \ref{assumption}).
\end{enumerate}
As a consequence, the set of candidate facilities of a $\dlt$-cover is in one-to-one correspondence to the edges and nodes of the network.
\end{observation}

With the above observations, we can already give a high-level description of the covering characterization behind our model. 
Namely, if we fix a set of facilities on some network edges, there would be some edges completely covered regardless of the exact facility locations within their edges. Some other edges would be partially covered from the left-end node and/or from the right-end node, and how much
depends on the actual facility locations. So we have variables that
specify the exact facility locations on each edge. Finally, we stipulate that the cover from the left and the cover from the right better exceed the edge length. 
One difficulty is that the corresponding covering function
is not linear 
However, we show that such a function is a piece-wise linear function, which can be modeled by a MILP. The aim of the remainder of this section is to give a mathematical  specification of the above characterization.

\subsection{Covering conditions}
We give a sufficient and necessary condition that the network is $\dlt$-covered by installed facilities.
 
 \begin{proposition}
 \label{thm.cover}
Let $\cP=\{p_i\}_{p_i\in C(N)}$ be a finite set  of points in $C(N)$. An edge $e=(v_a, v_b) \in E$ is $\dlt$-covered by $\cP$ if and only if either there exists $p \in \cP\cap e$ or
 \begin{equation}
	\label{eq.covercond}
    	\max\{\dlt -  \min_{p \in \cP} d(v_a,p),0\} +  \max \{\dlt - \min_{p \in \cP} d(v_b,p), 0 \} \ge l_e.
	\end{equation}
Moreover, the set $\cP$ is a $\dlt$-cover of $N$ if and only if for each $e\in E$, either there exists $p \in \cP\cap e$ or \eqref{eq.covercond} is satisfied.
\end{proposition}
\begin{proof}  
If there exists $p\in \cP\cap e$, then $e$ is $\dlt$-covered by $p$ due to Assumption \ref{assumption}.
Otherwise, for each $i\in\{a,b\}$, let us consider $p^*_i \in \cP$ such that $d(v_i,p^*_i)=\min_{p \in \cP} d(v_i,p)$ and let  $\pi^*_i\in \Pi^*(v_i,p_i^*)$ be a shortest path between $v_i$ and $p_i^*$, i.e. $l_{\pi^*_i}= d(v_i,p^*_i)$. Condition \eqref{eq.covercond} can be rewritten as follows:
$$ \max\{\dlt - l_{\pi^*_a},0\} +  \max \{\dlt - l_{\pi^*_b}, 0 \} \ge l_e.$$
Note that $\max\{\dlt - l_{\pi^*_i},0\} $ represents the maximum length that can be $\dlt$-covered by $\cP$ (specifically, from $p^*_i$) after passing through $v_i$. Since the path(s) that $\dlt$-cover $e$ must contain $v_a$ and/or $v_b$, the edge is covered if and only if these ``maximum lengths'' for $v_a$ and $v_b$ add up to more than $l_e$.

\end{proof}
\begin{figure}
	\begin{subfigure}{0.33\textwidth}
	\centering
   	\includegraphics[scale=1]{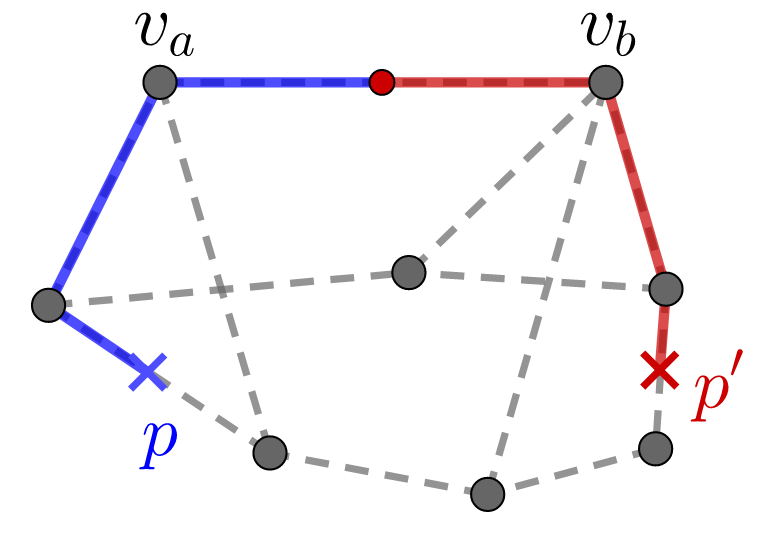}
   	\caption{By two points}
   	\label{fig:teoex-a}
	\end{subfigure}
	\begin{subfigure}{0.33\textwidth}
	\centering
   	\includegraphics[scale=1]{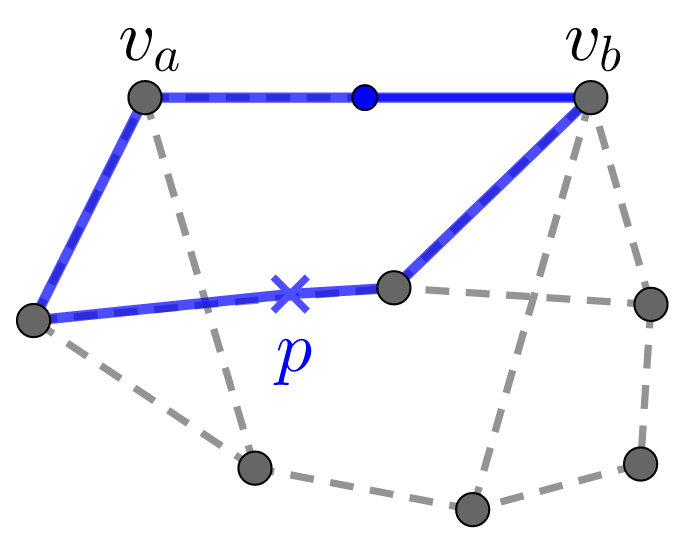}
   	\caption{By one point, through both ends}
   	\label{fig:teoex-b}
	\end{subfigure}
 	\begin{subfigure}{0.33\textwidth}
	\centering
   	\includegraphics[scale=1]{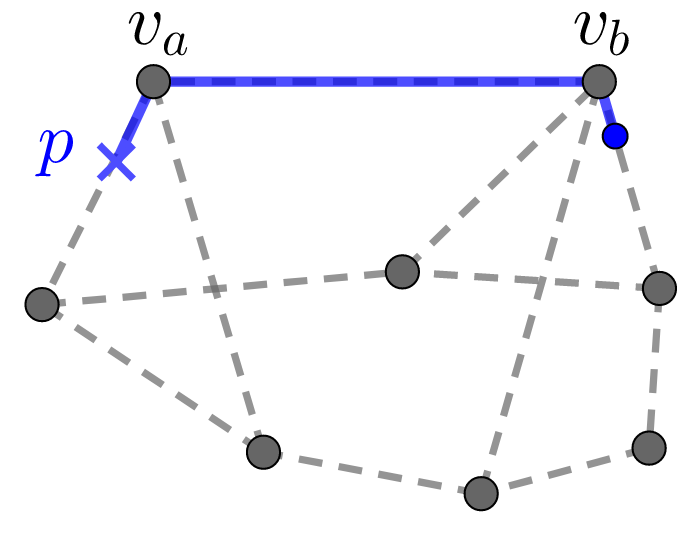}
   	\caption{By one point, through one end}
   	\label{fig:teoex-c}
	\end{subfigure}
	\caption{Covering of an edge $e=(v_a,v_b)\in E$}
	\label{fig:teoex}
\end{figure}
Figure \ref{fig:teoex} illustrates Proposition \ref{thm.cover}. It shows three ways of covering the same edge $e=(v_a,v_b)\in E$ for a given network. The edges of the network are depicted with dashed lines, while the different paths through which $e$ is covered are delimited with continuous bold traces. Facility locations are marked with the symbol `x'. Figure \ref{fig:teoex-a} depicts two facilities located at $p,p'$ that cover two portions of the edge, which contain $v_a$ and $v_b$ respectively. In this case, the $\min$ functions inside \eqref{eq.covercond} are attained respectively at $p$ and $p'$. In the middle,  Figure \ref{fig:teoex-b} shows a single location $p$ that covers $e$ through two different paths, which traverse $v_a$ and $v_b$ respectively. These paths form a cycle that contains $p$ and $e$. In this case,  the two $\min$ operations inside \eqref{eq.covercond} are attained at the same point, $p$. Finally, Figure \ref{fig:teoex-c} illustrates the case in which a single facility located at $p$ covers $e$ through one of its end nodes, $v_a$. Here, one of the $\max$ operators in \eqref{eq.covercond} is equal to zero ($p$ is further from $v_b$ than $\dlt$).

\subsection{Covering delimitation and simplification}
\label{sec.subc}

For an optimization or search problem, \textit{delimitation} refers to the reduction of the candidate space. As we have seen in \Cref{sub.disc}, FDS is studied in related works as a way for reduction of \problem\ (under some assumption), and hence it is a kind of delimitation. Instead of FDS, we consider a different delimitation that can be used for the general \problem, and which allows us to obtain a reduced MILP formulation.

The characterization in Proposition \ref{thm.cover} is based on the individual covering of every edge in the network. When considering possible locations to cover a fixed edge, we can restrict ourselves to its surroundings within the radius $\dlt$. Delimiting those parts of the network that could ``contribute'' to covering a particular edge or node reduces the search space. We  then introduce three kinds of delimitation: \textit{potential covers}, \textit{complete covers} and \textit{partial covers}.   We will represent the covering condition under such delimitation. In effect, the covering condition has a simplified form compared to its general form in \Cref{thm.cover}, and an adequate preprocessing procedure can reduce and strengthen our MILP model.

We find that, for every node, there may exist potential covers, i.e., a set of edges and nodes where, if a facility is located, it can \textit{possibly} $\dlt$-cover this node.  The potential covers in the following definition delimit the edges and nodes of $C(N)$ that can contribute to covering a particular node of the network.
\begin{definition}\label{defa.linksetsv}
For each $v \in V$, the potential covers of $v$ are the candidate facility locations to cover $v$:
\begin{align}
 \cE(v)&:= \{e'=(v'_a,v'_b) \in E:\:   d(v,v'_{i}) \le \dlt \textup{ for some } i\in \{a,b\} \} \nonumber\\
 \cV(v)&:= \{v' \in V: \:  d(v,v') \le \dlt\} \nonumber\\
 \cF(v)&:= \cE(v) \cup \cV(v). \nonumber
\end{align}
\end{definition}
Clearly, $v$ is not reachable within the radius $\dlt$ for any facility installed outside $\cF(v)$. Regarding the covering of edges, an edge incident to $v$ could be covered by some of the facilities in $\cF(v)$.

We find that, for every edge $e$, there may exist complete covers, i.e., a set of edges and nodes where, if a facility is located, it can \textit{always} $\dlt$-cover the edge $e$, regardless of the exact facility location. Once there is a facility within a complete cover, the whole edge $e$ is guaranteed to be covered by this facility. The following definition serves to delimit the network to such complete covers.

\begin{definition}
\label{def.linksetse}
For each $e=(v_{a}, v_{b}) \in E$, the  complete covers of $e$ are the candidate facility locations that can completely cover $e$:
\begin{align}
\cE_{\mathrm{c}}(e)& :=\{e' \in E:   \forall p' \in e', \forall  p \in e,\:   d(p,p') \le \dlt\}\nonumber\\
   \cV_{\mathrm{c}}(e)&:= \{v' \in V:   \forall  p \in e,\: d(p,v') \le \dlt \}\nonumber\\
  \cF_{\mathrm{c}}(e) &:= \cE_{\mathrm{c}}(e) \cup \cV_{\mathrm{c}}(e).
\end{align}
\end{definition}
If a facility is placed at $\cF_{\mathrm{c}}(e)$ (either at a node in $\cV_{\mathrm{c}}(e)$ or at a point on an edge belonging to $\cE_{\mathrm{c}}(e)$), we can immediately conclude that $e$ is $\dlt$-covered.
Note that any facility placed at $e'$ can completely cover $e$ if and only if any facility placed at $e$ can completely cover $e'$. That is,  $\cE_{\mathrm{c}}$ is symmetric over $E$. On the other hand, it is obvious that $e\in \cE_{\mathrm{c}}(e)$,  and $v_{a} ,v_{b}\in \cV_{\mathrm{c}}(e)$, for all $e=(v_a,v_b)\in E$.

Given a node $v$, we can characterize the complete covers of the incident edges to the node $v$. This helps us refine the potential covers of this node.  The following definition  identifies those candidate facility locations in the potential covers of the node $v$ that cannot completely cover any incident edges to $v$.
\begin{definition} \label{defa.linksetsvincident}
We define the following sets for each $v\in V$:
\begin{align}
\cE_{\mathrm{p}}(v)&:=\{e'\in \cE(v):\: \exists e\in E(v),\: e'\notin\cE_{\mathrm{c}}(e)\} \nonumber\\  \cV_{\mathrm{p}}(v)&:=\{v'\in \cV(v):\: \exists e\in E(v),\: v'\notin\cV_{\mathrm{c}}(e)\}   \nonumber\\
\cF_{\mathrm{p}}(v)&:=\cV_{\mathrm{p}}(v)\cup\cE_{\mathrm{p}}(v).
\end{align}
We call these sets the \emph{partial covers} of $E(v)$, where $E(v):=\{e\in E:\: v\in e\}$ is the set of incident edges to $v$.
\end{definition}
The set $\cF_{\mathrm{p}}(v)$ contains those candidate locations that can contribute to partially covering some of the edges in $E(v)$. Note that, if a facility is placed at $\cF(v)\setminus \cF_{\mathrm{p}}(v)$, then this facility completely covers  $E(v)$.

Definitions  \ref{defa.linksetsv}, \ref{def.linksetse}, and \ref{defa.linksetsvincident}  provide us with a refined covering condition. Indeed, the following proposition is a consequence of Proposition \ref{thm.cover} and the aforementioned definitions, and will be used to characterize coverings in the MILP formulation presented in Section \ref{sec:milp}.

\begin{proposition}
\label{cora.covercond}
A finite set $\cP$ of points in $C(N)$  is a $\dlt$-cover of $N$ if and only if,
for each $e=(v_a,v_b) \in E$, either $\cP \cap \cF_{\mathrm{c}}(e) \ne \varnothing$ or
\begin{equation}
\label{eqa.coverdistfull}
	\sum_{i \in \{a,b\}} \max \big\{0,\dlt -  \min_{p \in   \cP \cap \cF_{\mathrm{p}}(v_i)} d( v_{i}, p)\big \} \ge l_e.
\end{equation}
Moreover,
\begin{equation}\nonumber
\begin{split}
\min_{p \in \cP \cap \cF_{\mathrm{p}}(v_i)}d(v_{i}, p)
   = &\min \left\{\min_{v' \in \cP \cap \cV_{\mathrm{p}}(v_i)} d(v_{i}, v'), \min_{p \in  \cP \cap (\cE_{\mathrm{p}}(v_i)\setminus  V)}  d(v_{i}, p) \right\}, \textup{ for } i=a,b.
\end{split}
\end{equation}
\end{proposition}
\begin{proof}
For $i \in \{a,b\}$, the equality $\cF(v_i)=\cF_{\mathrm{p}}(v_i)\cup\cF_{\mathrm{c}}(e)$ holds by definition, which gives the new necessary and sufficient covering condition.

For the second statement of the proposition, we have $\cF_{\mathrm{p}}(v_i)=\cV_{\mathrm{p}}(v_i)\cup \cE_{\mathrm{p}}(v_i)$ from Definition \ref{defa.linksetsvincident}. Then, it suffices to see that we can take $p \in  \cP \cap (\cE_{\mathrm{p}}(v_i)\setminus  V)$ in the second inner $\min$ operator of the right-hand side instead of $p \in  \cP \cap \cE_{\mathrm{p}}(v_i)$. Let $e'=(v'_a,v'_b)\in  \cE_{\mathrm{p}}(v_i)$. We prove that the end nodes of $e'$ can be excluded from the second inner $\min$ operator.  Let $i' \in \{a,b\}$, we denote $\bar{i'}=b$ if $i'=a$ and $\bar{i'}=a$ if $i'=b$.  First, by definition, $d(v_i,v'_{i'})\leq \dlt$ for some $i'\in\{a,b\}$. Therefore, $v'_{i'}\in \cV_{\mathrm{p}}(v_i)$, and we can exclude it from the second inner $\min$ operator (this node is already considered by the first inner $\min$ operator). We consider now the other end node of $e'$. If $d(v_i,v'_{\bar i'})\leq \dlt$ then, similarly, $v'_{\bar{i}'}$ can be excluded from the second inner $\min$. Otherwise, we know that the outer $\min$ is not attained at $v'_{\bar i'}$, as  $d(v_i,v'_{i'})\leq \dlt < d(v_i,v'_{\bar i'})$, thus $v'_{\bar i'}$ can be disregarded.
\end{proof}

\begin{remark}
The covering conditions described both in Propositions \ref{thm.cover} and \ref{cora.covercond} would be also applicable if only a subset of $E$, $J\subseteq E$, is to be covered. Indeed, these covering conditions are based on the individual coverage of the edges, so it would be sufficient to apply them just to the edges in $J$. As a consequence, our methods, including the MILP formulation and algorithms presented in the next sections, apply to this more general version of the \problem.
\end{remark}

In order to exploit the newly defined potential, complete, and partial covers in our formulation, from a practical viewpoint, we need to have some characterizations that can operate in a computer. Definition \ref{defa.linksetsv}, which introduces potential covers, satisfies this requirement. Indeed, it just depends on distances between pairs of nodes, which we can easily calculate. Conversely, Definition \ref{def.linksetse} presents complete covers with a condition that must be satisfied by ``the infinitely many points of an edge'', which is not directly computable. Finally, the elements in the partial covers defined by Definition \ref{defa.linksetsvincident} can be easily calculated once both potential and complete covers are known.

\subsection{Characterization of complete covers}

In the following, we focus on characterizing the complete covers, which will be useful for our MILP formulation and tailored algorithms, (see forthcoming Sections \ref{sec:milp} and \ref{sec:prepro}). To begin with, we note that Proposition \ref{thm.cover} already gives us a characterization of the nodes in $\cV_{\mathrm{c}}(e)$. Indeed, it is easy to observe that, for a given $e\in E$, $v\in \cV_{\mathrm{c}}(e)$ if and only if $\cP:=\{v\}$ $\dlt$-covers $e$. We then focus on the sets $\cE_{\mathrm{c}}(e)$. On the one hand, it is clear that, for every edge $e=(v_a,v_b)\in E$,
$$\cE_{\mathrm{c}}(e)\subseteq \cE(v_a)\cap \cE(v_b).$$
Moreover, if we define $\cE_{\mathrm{c}}(v)\subseteq \cE(v)$ as follows, we have a tighter set containing $\cE_{\mathrm{c}}(e)$.
\begin{definition}\label{def.cEc}
The edges that can completely cover a node $v \in V$ are:
$$\cE_{\mathrm{c}}(v) :=\{e' \in E:   \forall p' \in e',\:   d(v,p') \le \dlt\}.$$
\end{definition}

It is clear that the following observation holds.
\begin{observation}\label{obser}
For any $e=(v_a,v_b)\in E$, $\cE_{\mathrm{c}}(e)\subseteq \cE_{\mathrm{c}}(v_a)\cap \cE_{\mathrm{c}}(v_b)$.
\end{observation}
With this observation, we can limit the search of $\cE_{\mathrm{c}}(e)$ in the set $\cE_{\mathrm{c}}(v_a)\cap \cE_{\mathrm{c}}(v_b)$, as we will show in \Cref{sec:prepro}, the latter set is easy to compute.

We recall that Definition \ref{def.cEc} is somewhat the inverse of Definition \ref{def.linksetse}. That is,  $e'\in \cE_{\mathrm{c}}(v)$ if and only if $v \in \cV_{\mathrm{c}}(e')$. We present a set of intermediate statements in Definition \ref{def.function}, Lemma \ref{lem.covercE}, and Lemma  \ref{lem.pieces}, which allow us to describe the edges in the complete cover set $\cE_{\mathrm{c}}(e)$ as the main result in Proposition \ref{prop.chara.cEe}.

\begin{definition}\label{def.function}
Let $v\in V$ be a node and $e'=(v'_a,v'_b)\in E $ be an edge. For all $q \in [0, l_{e'}]$, we define the following functions:
\begin{align}
	d_v(q) &:= \min \{d(v,v'_a)+q, d(v,v'_b)+l_{e'}-q\} \nonumber\\
	r_v(q) &:= 
	\max\{\dlt -  d_v(q),0\}, \nonumber
\end{align}
and constant:
$$Q_{ve'}:=(d(v,v'_b)+l_{e'} - d(v,v'_a))/2,$$  which satisfies the following equation:
$$d(v,v'_a)+Q_{ve'}=  d(v,v'_b)+l_{e'}-Q_{ve'}.$$
\end{definition}
The function $d_v(q)$ represents the distance between $v$ and a point $p'\in e'$ such that $q=l(v'_a,p')$, where $l(v'_a,p')$ measures the length of the continuum $(v'_a,p')\subseteq e'$.
The inner terms in the minimization that defines $ d_v(q)$ coincide for $q=Q_{ve'}$. Informally, $Q_{ve'}$ is the ``bottleneck" coordinate on $e'$ (see original definition in \cite{obnoxious}), for which the distance to $v$ is the same if we go through $v'_a$ or $v'_b$. Indeed, since $|d(v,v'_b) - d(v,v'_a)| \le l_{e'}$, it follows that $0 \le Q_{ve'} \le l_{e'}$.
Note that $d_v(q)$,  $r_v(q)$, and $Q_{ve'}$ depend also on the edge $e'$. 

\begin{lemma} \label{lem.covercE}
Let $e=(v_a,v_b)$. An edge $e'$ is in $\cE_{\mathrm{c}}(e)$ if and only if $ r_{v_a}(q)+r_{v_b}(q)\geq l_e$ for all $q\in [0,l_{e'}]$.
\end{lemma}
\begin{proof}
$e' \in \cE_{\mathrm{c}}(e)$ if and only $e$ is $\dlt$-covered by any point $p' \in e'$.
Take $\cP = \{p'\}$ in \Cref{thm.cover}, $e$ is $\dlt$-covered by $\cP$, if and only if
 \begin{equation*}
 \max\{\dlt -  d(v_a,p'),0\} +  \max\{\dlt -  d(v_b,p'),0\} \geq l_e
\end{equation*}
holds for all $p'\in e'$. Let $q=l(v'_a,p'),\: q \in [0, l_{e'}]$, be the measure of the sub-edge $(v'_a,p')\subseteq e'$. Then, the observation that $d_{v_i}(q)=d(v_i,p')$ for $i=\{a,b\}$ completes the proof.
\end{proof}

We present the following lemma without proof. The lemma is a direct consequence of Definition \ref{def.function}.
\begin{lemma}\label{lem.pieces}
Let $v\in V$ and $e'=(v'_a,v'_b)\in E$. Then, $0 \le Q_{ve'} \le l_{e'}$, the function $d_v(q)$ is increasing when $q\in [0, Q_{ve'}]$, and it decreases for $q\in[Q_{ve'}, l_{e'}]$. Moreover, $d_v(q)$ admits the following piece-wise linear representation:
\begin{equation*}
	d_v(q) = \begin{cases}
	d(v,v'_a)+q & \textup{ if }q \le Q_{ve'}, \\
	d(v,v'_b)+l_{e'}-q & \textup{ if }q \ge Q_{ve'}. \\   
	\end{cases}
\end{equation*}
\end{lemma}

Note that with  \Cref{obser}, to find $\cE_{\mathrm{c}}(e)$, we can check \Cref{lem.covercE} for $e' \in \cE_{\mathrm{c}}(v_a) \cap \cE_{\mathrm{c}}(v_b)$: if $ r_{v_a}(q)+r_{v_b}(q)\geq l_e$, then $e'$ is in  $\cE_{\mathrm{c}}(e)$. Then, we consider $v \in V$, and we want to characterize $r_v$ for all points on edges $e'\in E$ such that $e' \in \cE_{\mathrm{c}}(v)$. We present the following result without proof.

\begin{lemma}
Given $v \in V$, and $e'=(v'_a,v'_b) \in \cE_{\mathrm{c}}(v)$. If $e'\in \cE_{\mathrm{c}}(v)$, then $r_v(q)$ admits the following piece-wise linear representation:
\begin{equation*}
	r_v(q) = \begin{cases}
	\dlt - (d(v,v'_a) + q) & \textup{ if }q \le Q_{ve'}, \\
	\dlt- (d(v,v'_b)+l_{e'}-q) & \textup{ if }q \ge \ Q_{ve'}. \\   
	\end{cases}
\end{equation*}
\end{lemma}

Finally, we have a tractable version of \Cref{lem.covercE}.

\begin{proposition}\label{prop.chara.cEe}
Let $e=(v_a,v_b)$. An edge $e'\in \cE_{\mathrm{c}}(v_a)\cap \cE_{\mathrm{c}}(v_b)$ is in $\cE_{\mathrm{c}}(e)$ if and only if
$$r(q)\geq l_e\quad \textup{for all}\quad q=Q_{v_ae'}, Q_{v_be'},$$
where $r(q):=r_{v_a}(q)+r_{v_b}(q)$.
\end{proposition}
\begin{proof}
From Lemma \ref{lem.covercE}, $e' \in \cE_{\mathrm{c}}(e)$ if and only if $\min_{q \in [0, l_{e'}]} r(q) \ge l_e$. Due to Lemma \ref{lem.pieces}, the minimum argument must be some of the breakpoints in the piece-wise linear description of $r_v(q)$. Then, it suffices to check
\begin{equation*}
	\min \{r(q):\: q \in \{0, Q_{v_ae'}, Q_{v_be'}, l_{e'}\}\} \ge l_e.
\end{equation*}

Since $e'\in \cE_{\mathrm{c}}(v_a)\cap \cE_{\mathrm{c}}(v_b)$, $r_{v_i}(0)\geq l_e$ and  $r_{v_i}(l_{e'})\geq l_e$  always holds for all $i\in\{a,b\}$.
 It suffices thus to check  the following condition
\begin{equation*}
	\min \{r(q):\: q \in \{ Q_{v_ae'}, Q_{v_be'} \}\}  \ge l_e.
\end{equation*}
Then, the result follows.
\end{proof}

\section{MILP formulation}\label{sec:milp}

We present the variables of the MILP formulation first. For each $v \in V$, there is one candidate facility (fixed location); and for each $e \in E$, there is another one (location within the interior of $e$, $e\setminus\{v_a,v_b\}$). Then, the finite set $\cF \coloneqq E \cup V$ will be used to index the candidate facilities.
In our MILP formulation, there are two decisions associated with each  $f\in \cF$. One is to decide if a facility is installed at $f$. The second is only necessary for those facilities installed at the interior of edges, and consists in determining their locations within the corresponding edges.

To represent the first of the above decisions, we define the following binary variables, which we call the \emph{placement variables}:
$$y_f=1 \textup{ if a facility  is installed at } f,\quad  \textup{ for all }f \in \cF.$$
We identify the set of installed facilities with $\cF_1 = \{f \in \cF: y_f = 1\}$. To represent the second of the above decisions, we define the following continuous variables, which we name the \emph{coordinate variables}:
$$q_e=
\begin{cases}
l(v_a,p) &\textup{ if } y_e=1 \textup{ and a facility is installed at } p\in e\\
0 &\textup{ otherwise}
\end{cases}
\qquad \textup{ for all }e=(v_a,v_b)\in E.$$

We use $v'$ and $e'$ to denote nodes and edges where facilities are installed,   and use $v$ and $e$ to denote  nodes and edges to be covered, respectively. We refer to $v_a, v_b$ as the end nodes of $e=(v_a,v_b)\in E$; given $e'=(v'_a,v'_b) \in E$, we refer to $v'_a, v'_b$ as its end nodes.

We use the necessary and sufficient condition of $\dlt$-covering in Proposition \ref{cora.covercond}, and the second result in this proposition regarding the distance function. Other than the placement and coordinate variables, some additional variables are used, which we present next.

\setlength{\tabcolsep}{3pt}
\begin{tabular}{lll}
$w_e\in\{0,1\}$ & $=1$ if $\cF_{\mathrm{c}}(e)\cap \cF_1\neq \varnothing$  & $ e \in E$;\\
$r_{v}\in[0, +\infty)$ &
& $ v \in V$; \\
$x_{v}\in\{0,1\}$ & $=1$ if,  for all $e\in E(v)$, $ \cF_{\mathrm{c}}(e) \cap \cF_1 \ne \varnothing$ & $ v \in V$;\\
$z_{vv'}\in\{0,1\}$ & $=1$ if, for all $e=(u,v)\in E(v)\st w_e=0$, $\max\{0,\dlt-d(v,v')\}+r_u\geq l_e$ & $ v \in V$, $v' \in \cV_{\mathrm{p}}(v)$;\\
$z_{ve'{i'}}\in\{0,1\}$ & $=1$ if, for all $e=(u,v)\in E(v)\st w_e=0$, $\max\{0,\dlt- \tau_{ve'i'}(q_{e'})\}+r_u\geq l_e$ & $ v \in V$, $(e',i') \in \cEI_{\mathrm{p}}(v)$,
\end{tabular}

where $$ \tau_{ve'i'}(q) := d(v,v_{i'})+ \mathbb{1}_{i'=a}q +\mathbb{1}_{i'=b}(l_{e'} - q),\quad \textup{and}\quad \cEI_{\mathrm{p}}(v) := \{(e'=(v'_a,v'_b),i') \in \cE_{\mathrm{p}}(v) \times \{a,b\}:\: d(v,v'_{i'})\leq \dlt\}.$$ We sometimes refer to $r_v$ as the ``residual cover'' at node $v$, since it represents the maximum remaining length that can be covered after reaching $v$ from ``a sufficiently close'' facility. Remembering \Cref{def.function}, these variables must satisfy:
$$r_v\leq r_v(q_e), \: \textup{ for all } e \in \cF_1 \cap (E \smallsetminus E(v)).$$
In an optimal solution, it can be $r_v=r_v(q_e)$ for an edge $e\in E$ as stated above. In this case, $r_v$ is the maximum remaining length that can be covered after reaching $v$ from the closest facility. However, we do not impose this equality in our formulation, since it is enough for guaranteeing the coverage of $e=(v_a,v_b)\in E$ that the sum of the residuals $r_{v_a}+r_{v_b}$ exceeds $l_e$ (see Proposition \ref{cora.covercond}). That is, if $l_e$ is already exceeded by $r_{v_a}+r_{v_b}$  for some $r_{v_i} < r_{v_i}(q_e)$, then the coverage condition of Proposition \ref{cora.covercond} will hold.

We denote by $M_{*}$ the sufficiently large big-M constant associated with index $*$, the value of which will be determined later. Our formulation of the \problem\ reads as follows:
\begin{subequations}
	\label{milpa}
	\begin{align}
   	\min & \sum_{f \in \cF} y_f \label{milpa.obj} \\
	\st   & w_{e} \ge y_{f}   &  e \in E, f \in \cF_{\mathrm{c}}(e) \label{milpa.completecover1}\\
      	& w_{e} \le  \sum_{f \in \cF_{\mathrm{c}}(e)}y_{f} 	&  e \in E \label{milpa.completecover2}\\
      	& x_{v}\geq 1-\sum_{e\in E(v)} (1-w_e) &  v\in V \label{milpa.enforcex1}\\
      	& x_{v}\leq w_e &  v\in V, e\in E(v) \label{milpa.enforcex2}\\
	&  y_{v'_{i'}} + y_{e'} \le 1 &   e' \in E, i' \in \{a,b\} \label{milpa.optimale}\\
	& q_{e'}\leq l_{e'}y_{e'} & e'\in E  \label{milpa.coord} \\
	&  l_e (1 - w_e) \le  r_{v_a} + r_{v_b}	& e \in E \label{milpa.covere}\\
  &  x_{v} +\sum_{v' \in \cV_{\mathrm{p}}(v)}  z_{vv'}+ \sum_{(e', i') \in \cEI_{\mathrm{p}}(v)} z_{ve'i'} = 1  &   v\in V \label{milpa.lcoversos} \\
  &   z_{vv'} \le y_{v'}  &  v \in V,  v' \in \cV_{\mathrm{p}}(v) \label{milpa.indicatorv}\\  
	&   z_{ve'i'} \le y_{e'}  &  v \in V,  (e',i') \in \cEI_{\mathrm{p}}(v) \label{milpa.indicatorei}\\  
 &   r_{v} \le M_v (1-x_{v}) &  v\in V \label{milpa.wbdl}\\
 &   r_{v} \le M_{vv'}(1 - z_{vv'})  +  \dlt - d(v,v') & v \in V ,  v' \in \cV_{\mathrm{p}}(v) \label{milpa.coverdist-node} \\
 &   r_{v} \le M_{ve'i'}(1 - z_{ve'{i'}})  +  \dlt - \tau_{ve'i'}(q_{e'}) & v \in V ,  (e',i') \in \cEI_{\mathrm{p}}(v) \label{milpa.coverdist-edge} \\
   &  y_f, w_e  \in \{0,1\} &   f \in \cF,  e \in E \label{milpa.varyw}\\
   & x_{v}, z_{vv'}, z_{ve'{i'}} \in \{0,1\}  & v \in V, v' \in \cV_{\mathrm{p}}(v),(e',i') \in \cEI_{\mathrm{p}}(v)\label{milpa.varindicator}\\
  	& q_{e'}, r_{v} \geq 0 &  e'  \in E, v \in V. \label{milpa.varcontinuous}
	\end{align}
	\end{subequations}
Constraints \eqref{milpa.completecover1} and \eqref{milpa.completecover2} model the logic or constraint $w_e = \lor_{f \in \cF_{\mathrm{c}}(e)} y_f$.  Constraints \eqref{milpa.enforcex1} and \eqref{milpa.enforcex2} enforce the logic constraint  $x_v = \land_{e \in E(v)} w_e$, that is, $x_v$ is the product of the $w_e$ variables such that $e\in E(v)$.  
Constraints \eqref{milpa.optimale} prevent two facilities in a solution from being installed respectively at the interior of an edge and one of their end nodes.
Constraints \eqref{milpa.coord} bound the coordinate variables with the corresponding edge length, and set them to zero if no facility is located at its interior.
The covering condition in Proposition \ref{cora.covercond} is enforced by \eqref{milpa.covere}. If $w_e=1$, then the condition is satisfied ($e$ is covered by $\cF_{\mathrm{c}}(e)$). Otherwise, the inequality \eqref{eqa.coverdistfull} of the proposition has to be satisfied. The rest of the constraints of the model \eqref{milpa.lcoversos}-\eqref{milpa.coverdist-edge}, together with variables $r$, $x$, $q$, and $z$, aim at modeling \eqref{eqa.coverdistfull}. To begin with, \eqref{milpa.lcoversos} impose that, for each $v\in V$,  one of the following statements holds:
\begin{enumerate}
\item[i)] All incident edges to $v$, $e\in E(v)$, are completely covered by facilities placed at their complete covers, $\cF_{\mathrm{c}}(e)$  ($w_e=1$ for all $e\in E(v)$, $x_{v}=1$).
\item[ii)] A sufficiently close facility to $v$ is installed at $v'\in \cV_{\mathrm{p}}(v)$ ($z_{vv'}=1$), that is, $$\max\{0,\dlt-d(v,v')\}+r_u\geq l_e \qquad \forall u\in V \st (u,v)=e\in E \textup{ and }w_e=0;$$
\item[iii)]  A sufficiently close facility to $v$ is installed at $e'\in\cE_{\mathrm{p}}(v)$ and $v$ is reached through $v'_{i'}$ of $e'$ ($z_{ve'i'}=1$), that is, $$\max\{0,\dlt-\tau_{ve'i'}(q_{e'})\}+r_u\geq l_e \qquad \forall u\in V \st (u,v)=e\in E \textup{ and }w_e=0.$$
\end{enumerate}
If the case i) above holds, then the covering condition in \eqref{milpa.covere} is satisfied for all $e\in E(v)$, regardless of the value of the residual cover variables. Otherwise, suppose that $x_v=0$ and $w_e=0$ for some $e\in E(v)$. In this case, the corresponding constraint \eqref{milpa.covere} is ``active'', that is,  the inequality \eqref{eqa.coverdistfull} of Proposition \ref{cora.covercond} has to be satisfied for $e$. Since $x_v=0$, constraints \eqref{milpa.lcoversos} impose that there is a facility among those installed at $\cF_{\mathrm{p}}(v)$ that is sufficiently close one to $v$. This facility is the one bounding the residual variables $r_v$ (see constraints \eqref{milpa.coverdist-node}-\eqref{milpa.coverdist-edge}), which represent the terms in the left-hand side of \eqref{eqa.coverdistfull}.
Constraints \eqref{milpa.indicatorv} (resp. \eqref{milpa.indicatorei}) ensure that $z_{vv'}$ (resp. $z_{ve'i'}$) can be one only if  facility  is installed at $v'$ (resp. $e'$). Due to \eqref{milpa.lcoversos}, for every fixed node $v\in V$, at most one of the constraints in \eqref{milpa.wbdl}-\eqref{milpa.coverdist-edge} will be active. If $x_{v}=1$, \eqref{milpa.wbdl} enforces $r_{v}=0$. Indeed, all the covering conditions \eqref{milpa.covere} are ``inactive'' and $r_v$ is not needed to guarantee the coverage of any $e\in E(v)$. Otherwise, if $x_{v}=0$,  \eqref{milpa.wbdl} reads $r_{v}\leq M_v$, where $M_v$ is a big-enough constant that does not restrict the value of the residual.
Finally, constraints \eqref{milpa.coverdist-node}-\eqref{milpa.coverdist-edge} bound $r_v$ by $\dlt-d(v,p)\geq 0$ for a sufficiently close facility to $v$ installed at $p$, when $x_v=0$. The constants $M_{vv'}$ and $M_{ve'i'}$ are assumed to be big enough so that the constraints in \eqref{milpa.coverdist-node}-\eqref{milpa.coverdist-edge} do not add anything to the model if $z_{vv'}$ or $z_{ve'i'}$ are zero, respectively. For instance, $M_v=M_{vv'}=\dlt$ and $M_{ve'i'}=\dlt+l_{e'}$ are valid values for these constants (we recall Assumption \ref{assumption}).  Section \ref{sec:strenght} presents refined values of these big-Ms.

We observe that the number of variables and constraints in \eqref{milpa} can be reduced. Namely, for each $v\in V$ and $e'=(v'_a,v'_b)\in\cE_{\mathrm{p}}(v)$, if $d(v,v'_a)+l_{e'}\leq d(v,v'_b)$ then $d(v,p)=d(v, v'_{a})+l(v'_a,p)$ for every $p\in e'$. Similarly, if $d(v,v'_b)+l_{e'}\leq d(v,v'_a)$ then $d(v,p)=d(v, v'_{b})+l(v'_b,p)$ always holds for all $p\in e'$. For such nodes and candidate facilities, we do not need both variables, $z_{ve'a}$ and $z_{ve'b}$, and corresponding constraints in  \eqref{milpa.coverdist-edge} (we know beforehand that one of these constraints would never be active if a facility is located at $e'$). Therefore, $\cEI_{\mathrm{p}}(v)$ would only contain one of the pairs $(e',a)$ or $(e',b)$.

\subsection{Comparative insights with respect to an existing MILP}\label{sec:comp}
To the best of our knowledge, the only existing MILP for the \problem\ was proposed in \cite{Hamacher20}. The authors used a similar observation to ours with respect to optimal $\dlt$-covers. They noted that every edge contains at most one facility. Indeed, in their setting,  if two facilities are located at both end-nodes of an edge $e=(v_a,v_b)$ one of them is considered to be ``hosted'' by an adjacent edge, $e'\in E(v_a)\cup E(v_b)$ (by optimality, neither $v_a$ nor $v_b$ is a leaf). Their location variables are indexed then by $E$.

However, this approach has issues with symmetry, which leads to  redundant solutions. Indeed, many different solutions to the MILP model represent the same facility locations to  the \problem, since there exist many combinations of the edges ``hosting'' the facilities that are located at nodes. Let us consider an example: let $v \in V$ be a node, $e \in E(v)$ be an incident edge, and $p \in e$ be a point on $e$. Even though the distance $d(p,v)$ is equal to $ \epsilon $ for a very small $\epsilon > 0$, the point $p$ is still located at $e$. However, when $d(p,v) = 0$, the point $p$ is located at the node $v$ and, consequently, at every edge in $E(v)$. When a facility is located at a node, the discontinuity there leads to the question: which node or edge do we choose to represent this facility? In our MILP model, we have a specific node facility variable (i.e., $y_v$) which  prevents edge facility variables (i.e., $y_e, e \in E(v)$) from ``hosting" facilities at $v$, thanks to constraint \eqref{milpa.optimale}.

On the other hand, we consider model size in terms of the number of variables.  A second main difference between MILP \eqref{milpa} and the MILP in \cite{Hamacher20} is that the latter uses binary variables to identify the two edges containing the facilities that cover a given edge. On the one hand, this yields variables and constraints of $\bO{|E|^3}$. On the other hand, multiple equivalent solutions arise when the edge in question can be covered by a single facility, as the authors commented themselves. Finally, the covering constraints in both formulations actually correspond to the same characterization of $\dlt$-cover, but are modeled in a slightly different way. Namely, the authors of \cite{Hamacher20} defined the ``residual covers'' for each edge (where a facility might be placed) and node of the network. Interested readers might consult \cite{Hamacher20} and the MILP therein, which we do not reproduce here for the sake of concision.

Nonetheless, Table \ref{tab:milps} shows a comparative summary of the two formulations, based on the number of variables and constraints. This summary considers an upper bound on the size of MILP \eqref{milpa}. That is, we take $\cF_{\mathrm{c}}(e)=\cF$, $\cV_{\mathrm{p}}(v)=V$, and $\cE_{\mathrm{p}}(v)=E$ for all $v \in V$ and $e\in E$---however, this would never be the case, as the partial and complete covers are complementary. On the other hand, Table \ref{tab:milps}  considers the MILP in \cite{Hamacher20} with $J=E$ (the set of edges to be covered).

\begin{table}[]
	\centering
	\begin{tabular}{c|c|c|c}
     	& \multicolumn{2}{c|}{Variables} & Constraints\\[0.15cm]
     	\hline
     	& Binaries & Continuous \\[0.15cm]
     	\hline
 	MILP \eqref{milpa}& $|V|^2+2(|V||E|+|V|+|E|)$ & $|V|+|E|$ & $|E|^2+|V|^2+5|E||V|+7|E|+3|V|$ \\[0.15cm]
 	\hline
 	MILP in \cite{Hamacher20}&  $|E|^3+3|V||E|+|E|$ & $3|V||E|+|E|$ & $3|E|^3+8 |E||V|+|E|$
	\end{tabular}
	\vspace{0.35cm}
	\caption{Comparative summary on MILP formulations for the \problem}
	\label{tab:milps}
\end{table}

\section{Strengthening}\label{sec:strenght}
In this section, we analyze modifications of the MILP \eqref{milpa} that can yield a tighter linear relaxation of this formulation. Namely, we tighten our big-M constraints \eqref{milpa.wbdl}-\eqref{milpa.coverdist-edge} by devising small constants $M_v$, $M_{vv'}$, $M_{ve'i'}$, $\dlt_{vv'}$, and $\dlt_{ve'i'}$. We also present several families of valid inequalities. Valid inequalities define conditions that have to be satisfied by any feasible solution, and yield tighter linear programming relaxations, see e.g. \cite{wolsey}.

\subsection{Constants tightening}
\label{sec:strenght.bd}
From the MILP formulation, it is easy to yield the following observation. For $v \in V$, it suffices for a facility $f \in \cF_{\mathrm{p}}(v)$ to contribute to the residual cover $r_v$ at most $U_v := \max_{e \in {E}(v)} l_e$. Indeed, $r_v$ aims at ensuring that the inequality \eqref{eqa.coverdistfull} of Proposition \ref{cora.covercond} is satisfied for all $e\in E(v)$. We define
\begin{flalign*}
  	\dlt_{vv'}&:= \min\{U_v + d(v,v'), \dlt\}, & \textup{for } v' \in \cV_{\mathrm{p}}(v), \:\textup{ and }\qquad\qquad\:\:\\
   	\dlt_{ve'i'} &:= \min\{U_v + \max_{q \in [0, l_{e'}]} \tau_{ve'i'}(q) , \dlt\}
	=  \min\{U_v+ d(v,v'_{i'})+l_{e'}, \dlt \},& \textup{for } (e'=(v'_a,v'_b),i') \in \cEI_{\mathrm{p}}(v).
\end{flalign*}

Since $U_v$ is a valid upper bound for the residual cover variable $r_v$, the big-Ms in the constraints \eqref{milpa.coverdist-node} and \eqref{milpa.coverdist-edge} should guarantee that
\begin{equation*}
\begin{split}
	M_{vv'}+ \dlt_{vv'} - d(v,v') \ge\, U_v,\\
  M_{ve'i'}+ \min_{q \in l_{e'}}(\dlt_{ve'i'}-\tau_{ve'i'}(q)) \ge &\, U_v.
\end{split}
\end{equation*}

Taking the minimums of the above big-Ms, we can now tighten the big-M constants of the MILP \eqref{milpa} as follows:
\begin{equation}
\begin{split}
	M_{v}:=&\, U_v \\
	M_{vv'}:= &\, U_v - (\dlt_{vv'} - d(v,v')) = \max \{0, U_v+d(v,v')-\dlt\}\\
  M_{ve'i'}:= &\, U_v - \min_{q \in l_{e'}}(\dlt_{ve'i'}-\tau_{ve'i'}(q))
   = U_v - \dlt_{ve'i'} + \max_{q \in l_{e'}} \tau_{ve'i'}(q) \\
  =& U_v - \dlt_{ve'i'} + d(v_i,v'_{i'})+l_{e'}
  =  \max \{0, U_v+d(v,v'_{i'})+l_{e'}-\dlt\}, \nonumber
\end{split}
\end{equation}
where the last equations in the definition of $M_{vv'}$ and $M_{ve'i'}$ follow from the definition of $\dlt_{vv'}$ and $\dlt_{ve'i'}$, respectively.

Consequently, the constraints \eqref{milpa.coverdist-node} and \eqref{milpa.coverdist-edge}  should be replaced by:
\begin{subequations}
	\begin{align}
    	&   r_{v} \le M_{vv'}(1 - z_{vv'})  +  \dlt_{vv'} - d(v,v') & v \in V ,  v' \in \cV_{\mathrm{p}}(v) \label{milpa.coverdist-node_ref} \\
 &   r_{v} \le M_{ve'i'}(1 - z_{ve'{i'}})  +  \dlt_{ve'i'} - \tau_{ve'i'}(q_{e'}) & v \in V ,  (e',i') \in \cEI_{\mathrm{p}}(v). \label{milpa.coverdist-edge_ref}
	\end{align}
\end{subequations}

\subsection{Valid inequalities}\label{sec:valid-ineq}

\subsubsection*{``Leafs'' inequalities}
If a node $v\in V$ has degree one, we can assume without loss of generality that no facility is located at $v$ nor at its incident edge. Indeed, an equivalent $\dlt$-cover could be built by just moving such a facility to the unique neighbor of $v$ in $N$. More than valid inequalities, the following are valid variable elimination:
\begin{eqnarray}
y_v=0;\: y_e=0 & \forall  v \in V \: \st deg(v)=1, \: e\in E(v).\label{const-valid2}
\end{eqnarray}

\subsubsection*{``Adjacent edges'' inequalities}
Consider a node $v\in V$ of degree two. If there is a facility at $v$, then no facility is placed at the edges incident to $v$ (we recall the model constraints \eqref{milpa.optimale}). Otherwise, we can assume that at most one facility is placed at these edges in an optimal solution, which can be enforced by the following valid inequalities:
\begin{eqnarray}
y_e+y_{e'}+y_{v}\leq 1 & \forall  e,e'\in E, \: e\neq e',\: s.t. \: e\cap e'=v \textup{ and } deg(v)=2.\label{valid-deg2}
\end{eqnarray}
\begin{figure}
	\begin{subfigure}{0.5\textwidth}
    	\centering
    	\includegraphics[scale=1.1]{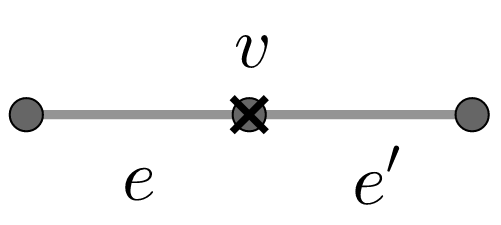}
    	\caption{A facility is located at $v$}
    	\label{fig:degree2a}
	\end{subfigure}
 	\begin{subfigure}{0.5\textwidth}
    	\centering
    	\includegraphics[scale=1.1]{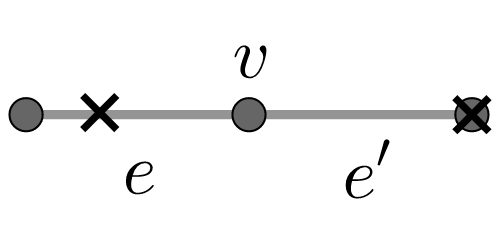}
    	\caption{Solution with two facilities}
    	\label{fig:degree2b}
	\end{subfigure}
	\caption{Illustration of valid inequalities \eqref{valid-deg2}}
	\label{fig:degree2}
\end{figure}
Figure \ref{fig:degree2} illustrates the above inequalities. Figure \ref{fig:degree2a} shows the case in which a facility is located at $v$. Otherwise, if two facilities are placed at $e$ and $e'$ respectively, we can build an equivalent solution by moving one of these facilities to the end node of the corresponding edge that is not $v$, as depicted in Figure \ref{fig:degree2b}. We recall that the last statement holds due to our Assumption \ref{assumption}.

\subsubsection*{``Neighborhood'' inequalities}
Let us now consider a node $v\in V$, and suppose that there are several facilities placed at different edges in $E(v)$ in a feasible solution. Take $e^*\in E(v)$ containing a facility $f^*$ such that $d(f^*,v)=\min \{d(f,v): f \textup{ is installed at }\:e\in E(v)\}$. The following proposition gives an equivalent feasible solution where the facilities at the edges $e\in E(v)$ such that $e\neq e^*$ are moved to the nodes.
\begin{proposition}\label{lemma:ineq}
Given a node $v \in V$ and an edge $e^* \in E(v)$, for any feasible solution $\hat{y}$ with several facilities placed at edges in $E(v)$, the following solution $y$ is feasible and $\sum_{f \in \cF} y_f \le \sum_{f \in \cF} \hat{y}_f$:
\begin{itemize}
	\item $y_{u}=1$ for all $u\in V$ such that $e=(u,v)\in E(v)$, $e \neq e^*$, and $\hat{y}_e=1$;
	\item $y_{e}=0$ for all $e\in E(v)$ such that $e \neq e^*$, and $\hat{y}_e=1$;
	\item $y_f=\hat{y}_f$ otherwise.
\end{itemize}
\end{proposition}
\begin{proof}
We denote by $N(v)$ the set of vertices adjacent to $v$. Consider the change of facilities from $\hat{y}$ to $y$. The facilities in the edges $E(v) \ \setminus \{e^*\}$ are `pushed' to the vertices $N(v)$. For $u \in N(v)$, if there already exists a facility at $u$, and there is another facility `pushed' to $u$, then these two facilities merge and they are accounted as one facility in  $y$. Hence, the number of facilities of the solution $y$ is at most that of the solution $\hat{y}$.

The proof then consists in showing that $y$ is feasible. We will show that all edges are covered. Let us consider $e\in E$. If $e$ was covered in $\hat{y}$ by facilities not placed at edges in $E(v)$ then it is still covered by these facilities in $y$. Suppose then that a facility placed at $e'\in E(v)$ with $e\neq e^*$ was covering $e$ (or part of $e$) in solution $\hat{y}$, and let $e'=(u,v)$. We distinguish two cases. First, if the facility at $e'$ was partially covering $e$ through node $u$, then it clearly covers at least the same part of $e$ in the new solution $y$ (where the facility is moved to $u$). Otherwise, suppose the facility at $e'$ was partially covering $e$ through node $v$. In this case, the facility at $e^*$ covers at least the same part of $e$ (it is closer to $v$). Since this facility remains unchanged in the new solution, we can guarantee that $e$ is still covered.
\end{proof}
\begin{figure}
	\begin{subfigure}{0.5\textwidth}
    	\centering
    	\includegraphics[scale=1.1]{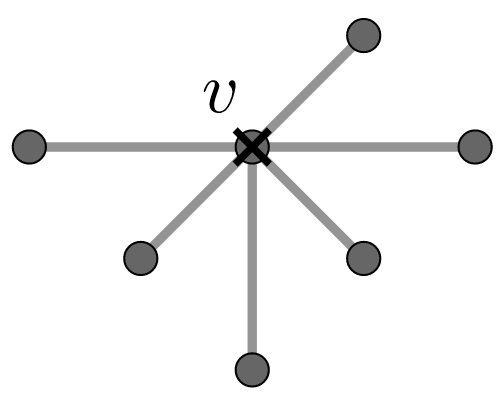}
    	\caption{A facility is located at $v$}
    	\label{fig:generala}
	\end{subfigure}
 	\begin{subfigure}{0.5\textwidth}
    	\centering
    	\includegraphics[scale=1.1]{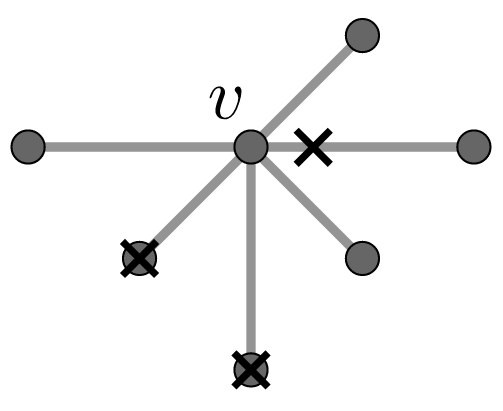}
    	\caption{Solution with several facilities}
    	\label{fig:generalb}
	\end{subfigure}
	\caption{Illustration of valid inequalities \eqref{ineq-neigh}}
	\label{fig:general}
\end{figure}
As a consequence of Proposition \ref{lemma:ineq}, the following inequalities are valid:
\begin{equation}
	\sum_{e\in E(v)} y_e \leq 1-y_v \qquad \forall v\in V. \label{ineq-neigh}
\end{equation}
Figure \ref{fig:general} illustrates the valid inequalities \eqref{ineq-neigh}. In particular, Figure \ref{fig:generalb} illustrates the equivalent solution given in Proposition \ref{lemma:ineq}. It is easy to observe that these new inequalities are a generalization of inequalities \eqref{valid-deg2}. Moreover, constraints \eqref{ineq-neigh} dominate the model constraints  \eqref{milpa.optimale}--- and are fewer.

\section{Network processing}\label{sec:prepro}
The network processing algorithm analyzes
the network $N$ to compute the parameters and sets needed to construct the MILP model \eqref{milpa}, which we recall next:

\begin{enumerate}
	\item $\cV_{\mathrm{c}}(e)$, $\cE_{\mathrm{c}}(e)$ for all edges $e \in E$;
	\item $\cV_{\mathrm{p}}(v), \cE_{\mathrm{p}}(v), \cEI_{\mathrm{p}}(v)$ for all nodes $v \in V$;
	\item $d(v,v')$ for all pairs of nodes $v,v' \in V$ such that $d(v,v') \le \dlt$.
\end{enumerate}

The above data is computed by Algorithms \ref{algo:vcover}, \ref{algo:mutual}  and \ref{algo:graph}. Algorithms \ref{algo:vcover} and \ref{algo:mutual} contain auxiliary functions, which are called within the main Algorithm \ref{algo:graph}. Algorithm \ref{algo:vcover} computes the sets  $\cE(v)$ and $\cV(v)$ (which are not directly used in the MILP but necessary to obtain $\cE_{\mathrm{p}}(v)$ and $\cV_{\mathrm{p}}(v)$), and the distances $d(v,v')$ for all $v,v'\in V$ such that $d(v,v')\leq \dlt$. Algorithm \ref{algo:vcover} also computes the sets $\cE_{\mathrm{c}}(v)$,
which will serve as intermediate sets to finally obtain $\cE_{\mathrm{c}}(e)$ in Algorithm \ref{algo:graph}.  The main task in Algorithm \ref{algo:graph} is to compute the sets  $\cV_{\mathrm{c}}(e)$ and $\cE_{\mathrm{c}}(e)$. To that aim, this algorithm calls both Algorithm \ref{algo:vcover} and the procedure ``mutual'' described in Algorithm \ref{algo:mutual}. Once  $\cV_{\mathrm{c}}(e)$ and $\cE_{\mathrm{c}}(e)$ are known, the computation of $\cV_{\mathrm{p}}(v)$ and $\cE_{\mathrm{p}}(v)$ in Algorithm \ref{algo:graph} easily follows by definition.

In the following, we present Algorithm \ref{algo:vcover}, which defines the function ``\ncalgo{}$(N,\dlt,s)$''. This function, for each source node $s \in V$, outputs: $\cE_{\mathrm{c}}(s)$, $\cE(s)$,  $\cV(s)$,
and $d(s,v)$ for all $v\in V$ such that $d(s,v)\leq \dlt$ (otherwise the algorithm outputs $d(s,v)=+\infty$). The algorithm  starts with empty sets $\cE_{\mathrm{c}}(s), \cE(s), U(s), \cV(s)$, where $U(s)$ is used for intermediate calculations. The set $Q$ denotes nodes whose shortest path (and distance) to $s$ are unknown, and it is initialized to $V$. In the course of the algorithm, $Q$ decreases, while $\cV(s)$ increases. In Lines \ref{algo:vcover.init.dist}-\ref{algo:vcover.start.covers}, the distance $d(s,v)$ and predecessor values $\prev_s(v)$ are initialized, for all $v\in V$. The while loop is an adaptation of the classic Dijkstra algorithm. Line \ref{algo:vcover.remove} selects the node $u$ with the shortest distance to $s$ among all unprocessed nodes, and removes it from $Q$.  If $d(s,u) > \dlt$, then none of the remaining nodes in $Q$ are reachable from $s$, and the search is pruned. Otherwise, the neighbors of $u$ that are still in $Q$ are inspected. For each $v\in Q\cap E(u)$, the edge $(u,v)$ is first added to $\cE(s)$. Then, the algorithm computes the length $\ell$ of a path from $s$ to $v$ that traverses $u$. If $\ell < d(s,v)$, then the distance and the predecessor for node $v$ are updated in Lines \ref{algo:vcover.updatedst}- \ref{algo:vcover.updatepred}. In addition, if $\ell < \dlt$,  node $v$ and edge $(u,v)$ are added to $\cE_{\mathrm{c}}(s)$ and $\cV(s)$ in Lines \ref{algo:vcover.ccovere} and \ref{algo:vcover.coverv}, respectively.  Otherwise, the edge $e$ is added to the undetermined set $U(s)$. Whether this edge belongs or not to the complete cover set $\cE_{\mathrm{c}}(s)$ is decided later on in the algorithm. Namely,  edges $e=(v_a, v_b)\in U(s)$ are processed in Lines \ref{algo:end1}-\ref{algo:end2}:  if $e$ can be jointly $\dlt$-covered  by $s$ from two sides, then $e$ is added to the complete cover $\cE_{\mathrm{c}}(s)$.

Algorithm \ref{algo:mutual} describes the procedure ``\malgo'', which determines, given $e=(v_a,v_b)\in E$ and a candidate edge for the complete cover $e'\in \cE_{\mathrm{c}}(v_a)\cap\cE_{\mathrm{c}}(v_b)$, whether $e' \in \cE_{\mathrm{c}}(e)$. This algorithm is based on Proposition \ref{prop.chara.cEe} in Section \ref{sec:prelim}.

Network processing Algorithm \ref{algo:graph} computes all the sets that are needed by the MILP formulation. The algorithm  starts with empty sets $\cE_{\mathrm{c}}(e), \cV_{\mathrm{c}}(e), \cE_{\mathrm{p}}(v), \cV_{\mathrm{p}}(v)$, for $e \in E$ and $v \in V$. In Line \ref{algo:graph.vcomplete}, the algorithm loops through all nodes  $v\in V$ and computes the function ``\ncalgo{}$(N,\dlt,v)$'', storing its output. Then, the algorithm calculates the sets $\cV_{\mathrm{c}}(e)$ for $e\in E$, by applying the symmetric relation between these sets and the sets $\cE_{\mathrm{c}}(v)$ from ``\ncalgo{}$(N,\dlt,v)$''. After that, in Line \ref{algo:graph.ecomplete}, the algorithm loops through all edges $e=(v_a,v_b) \in E$. It checks whether there is an edge $e' \in \cE_{\mathrm{c}}(v_a) \cap \cE_{\mathrm{c}}(v_b)$ such that $e' \in \cE_{\mathrm{c}}(e)$  (equivalently, $e \in \cE_{\mathrm{c}}(e')$) by calling the procedure "\malgo". Since $e' \in \cE_{\mathrm{c}}(e)$  if and only if $e \in \cE_{\mathrm{c}}(e')$, the loop only runs over pairs such that $e<e'$ (we assume a total order on the elements of $E$).  The loop starting in line \ref{algo:graph.vpartial}, iterates on each node $v \in V$ and looks for $v' \in \cV(v)$ such that there exists an $e \in E(v)$ but $e \notin \cE_{\mathrm{c}}(v')$. The nodes $v'$ found are added to $\cV_{\mathrm{p}}(v)$. Finally, the loop in line \ref{algo:graph.epartial} also iterates on $v\in V$, and looks for $e'=(v'_a,v'_b) \in \cE(v)$ such that there exists $e \in E(v)$ but $e \notin \cE_{\mathrm{c}}(e')$. Each edge found is added to $\cE_{\mathrm{p}}(v)$, and, right after that, the set $\cEI_{\mathrm{p}}(v)$ may be updated after checking the dominance rule described at the end of Section \ref{sec:milp}. We have the following complexity result for Algorithm \ref{algo:graph}.

\begin{proposition}
Let $D$ be an upper bound on the degree of the nodes of a connected network $N=(V,E)$.
The time complexity of the network processing Algorithm \ref{algo:graph} is  
${\bO{|E|^2+|V||E|(D+\log |V|)}}$.
\end{proposition}
\begin{proof}

We first analyze the time complexity of the procedure \ncalgo{} described in Algorithm \ref{algo:vcover}. The  main $\textbf{while}$ loop is a modification of the Dijkstra algorithm, and it can be implemented with time complexity ${\bO {(|E|+|V|)\log |V|}}$, see \cite{cormen2009introduction}. Therefore, the overall time complexity of \ncalgo{} is also ${\bO {(|E|+ |V|)\log |V|}}$.

The network processing Algorithm \ref{algo:graph} has four  outer  loops, and next we analyze the complexity of each outer loop.  The first outer loop runs the \ncalgo{} algorithm over the nodes, so its complexity is ${\bO {|V|(|E|+ |V|)\log |V|}}$. The second outer loop  runs the \malgo{} algorithm over the edges pairs, since the mutual algorithm has a constant time complexity, so the complexity of this loop is $\bO{|E|^2}$. The third outer loop is composed of three $\textbf{for}$ loops, and its time complexity is $\bO{D|V|^2}$, which has an upper bound $\bO{D|E||V|}$.  The last outer loop is composed of three $\textbf{for}$ loops, and its time complexity is $\bO{D|E||V|}$.  After summing up the complexity of these loops, we have that the total time complexity of Algorithm \ref{algo:graph} is  ${\bO{|V|((|E|+|V|)\log |V|)+|E|^2+|V||E|D}}$, or, equivalently,
${\bO{|V||E|\log |V|+|E|^2+|V||E|D}}$.

\end{proof}

\begin{algorithm}[htbp]
\SetAlgoLined
   \textbf{Input:} Network $N= (V, E, ||)$, cover range $\dlt > 0$,  a source $s \in V$\;
   \textbf{Output:} $\cE_{\mathrm{c}}(s)$, $\cE(s)$, $\cV(s)$, $d(s,v)$ for all $v\in V$ (returns $d(s,v)=+\infty$ if $d(s,v)>\dlt$)\;
   Initialize set $Q \gets V$\;
   Initialize sets $\cE_{\mathrm{c}}(s) \gets \varnothing$, $\cE(s) \gets \varnothing$, $U(s) \gets \varnothing$\;
   Initialize set $\cV(s) \gets \varnothing$\;
	\For{each node $v \in V$}{
    	$d(s,v) \gets +\infty$ \Comment*[r]{Unknown distance from $s$ to $v$} \label{algo:vcover.init.dist}
    	$\prev_{s}(v) \gets  \{\varnothing\}$  \Comment*[r]{Unknown predecessor of $v$} \label{algo:vcover.init.prev}
	}
	$d(s,s) \gets 0$\; 
	\label{algo:vcover.start.dists}
	add $s$ to $\cV(s)$\; 
	\label{algo:vcover.start.covers}
	\While{$Q$ is not empty}{
   $u \gets \argmin_{v \in Q} d(s,v)$\;
   remove $u$ from $Q$  \Comment*[r]{Take the closest node $u$ and remove it from $Q$} \label{algo:vcover.remove}
  \If{$d(s,u) > \dlt$} {
	$d(s,v) \gets +\infty$ for all $v\in Q$  \Comment*[r]{End of Dijkstra (all nodes in $Q$ are outside the covering radius)}\label{algo:vcover.break}
 	\Break
  }
   \For{each  $v\in Q$ s.t. $v\in E(u)$}{
   $e \gets (u,v)$\;
   add $e$ to $\cE(s)$ \Comment*[r]{Edge $e$ is in the potential cover set of $s$}
   $\ell \gets d(s,u) + l_e$ \Comment*[r]{Path from $s$ to $v$ that traverses $u$}\label{algo:vcover.distfound}
  	\If{ $\ell < d(s,v) $}{ \label{algo:vcover.term}
 	$d(s,v) \gets \ell$  \Comment*[r]{Update the distance to $v$}\label{algo:vcover.updatedst}
 	$\prev_{s}(v) \gets u$  \Comment*[r]{Update the predecessor of $v$}\label{algo:vcover.updatepred}
 	\uIf{ $\ell \le \dlt$}{
 	add $e$ to $\cE_{\mathrm{c}}(s)$ \Comment*[r]{Edge $e$ is in the complete cover set of $s$} \label{algo:vcover.ccovere}
	add $v$ to $\cV(s)$  \Comment*[r]{Node $v$ is in the potential cover set of $s$} \label{algo:vcover.coverv}
 	}
 	\Else{
 	add $e$ to $U(s)$ \Comment*[r]{Undetermined edge} \label{algo:vcover.undetermined}
 	}
   }
   }
   }
	\For{each edge $e = (v_a,v_b)$ in $U(s)$ \label{algo:end1} }{
  	\If{$v_a \in \cV(s)$ and $v_b \in \cV(s)$ and $ \dlt - d(s,v_a) + \dlt - d(s,v_b)  \ge l_e $}{
   	add $e$ to $\cE_{\mathrm{c}}(s)$  \Comment*[r]{Edge $e$ is completely covered} \label{algo:vcover.ccovere2}
  }   
 }\label{algo:end2}

 \caption{single node $\dlt$-cover algorithm: \ncalgo{}}
 \label{algo:vcover}
\end{algorithm}

\begin{algorithm}[htbp]
\SetAlgoLined
   \textbf{Input:} Edges $e=(v_a,v_b),e'=(v'_a,v'_b)$\ such that $e' \in \cE_{\mathrm{c}}(v_a)$ and $e \in \cE_{\mathrm{c}}(v_b)$.\\ 
   \textbf{Output:} Boolean value indicating whether $e' \in \cE_{\mathrm{c}}(e)$.\\\
\For{$i \in \{a,b\}$}{
 $Q_{v_ie'} \gets \frac{d(v_i,v'_b)+l_{e'} - d(v_i,v'_a)}{2}$\;  
}
 \For{$i \in \{a,b\}$}{
  \uIf{$q \le Q_{v_ie'}$ } {
 	$r_{v_i}(q) =  \dlt - (d({v_i},v'_a) + q)$ \;
  }
  \Else{
  	$r_{v_i}(q) =  \dlt- (d({v_i},v'_b)+l_{e'}-q)$ \;
  }
}
 \uIf{$r_{v_a}(Q_{v_ae'})+r_{v_b}(Q_{v_ae'})\ge e$ \textsc{and} $r_{v_a}(Q_{v_be'})+r_{v_b}(Q_{v_be'})\ge e$}{
 \Return \textsc{true}\;
 }
 \Else{
  \Return \textsc{false}\;
 }
 \caption{Edge mutual cover algorithm: \malgo}
 \label{algo:mutual}
\end{algorithm}

\begin{algorithm}[htbp]
\SetAlgoLined
   \textbf{Input:} Network $N= (V, E, ||)$\, and cover range $\dlt > 0$\;
   \textbf{Output:} $\cE_{\mathrm{c}}(e)$, $\cV_{\mathrm{c}}(e)$, $\cE_{\mathrm{p}}(v)$, $\cV_{\mathrm{p}}(v)$, $\cEI_{\mathrm{p}}(v)$, for all $e\in E$ and $v \in V$, and distance function $d$\;
 \For{ each node $v \in V$  \Comment*[r]{Computation of node complete covers $\cV_{\mathrm{c}}(e)$} \label{algo:graph.vcomplete}}{
   $\cE_{\mathrm{c}}(v), \cE(v), \cV(v),d(v,\cdot) \gets \ncalgo{(N, \dlt,v)}$ \;
 	\For{each edge $e \in \cE_{\mathrm{c}}(v)$}
   {
   add $v$ to $\cV_{\mathrm{c}}(e)$ \;
   }
 }
  \For{ each edge $e = (v_a, v_b) \in E$  \Comment*[r]{Computation of edge complete covers $\cE_{\mathrm{c}}(e)$} \label{algo:graph.ecomplete}}{
  	\For{ each edge $e' = (v'_a, v'_b) \in E$, $e<e'$, such that $e' \in \cE_{\mathrm{c}}(v_a)\cap \cE_{\mathrm{c}}(v_b) $}{
	\uIf{\malgo($e,e',d$)}{
	add $e'$ to $\cE_{\mathrm{c}}(e)$\;
 	add $e$ to $\cE_{\mathrm{c}}(e')$ \;
	}
	}
 }
 \For{each node $v \in V$ \Comment*[r]{Computation of node partial covers $\cV_{\mathrm{p}}(v)$} \label{algo:graph.vpartial}}{
   \For{each node $v' \in \cV(v)$}{
	\For{all $e \in E(v)$}{
	\uIf{$e \notin \cE_{\mathrm{c}}(v')$}{
 	add $v'$ to $\cV_{\mathrm{p}}(v)$\;
 	\Break
	}
	}
 }  
 }
 \For{each node $v \in V$ \Comment*[r]{Computation of edge partial covers $\cE_{\mathrm{p}}(v)$ and $\cEI_{\mathrm{p}}(v)$}  \label{algo:graph.epartial}}{
   \For{each edge $e'=(v'_a,v'_b) \in \cE(v)$}{
	\For{all $e \in E(v)$}{
	\uIf{$e' \notin \cE_{\mathrm{c}}(e)$
	}{
      	add $e'$ to $\cE_{\mathrm{p}}(v)$\;
 	\uIf{$d(v,v'_a)\leq \dlt$ \textsc{and}  $d(v,v'_a)\leq d(v,v'_b) + l_{e'}$}{add $(e',a)$ to $\cEI_{\mathrm{p}}(v)$\;}
 	\uIf{$d(v,v'_b)\leq \dlt$ \textsc{and}  $d(v,v'_b)\leq d(v,v'_a) + l_{e'}$}{add $(e',b)$ to $\cEI_{\mathrm{p}}(v)$\;}
 	\Break
	}
	}
 }  
 }
 \caption{Network processing algorithm}
 \label{algo:graph}
\end{algorithm}

\section{A reduced formulation for networks with long edges}\label{sec:reduced-milp}
\begin{figure}
	\begin{subfigure}{\textwidth}
    	\centering
    	\includegraphics[scale=1.15]{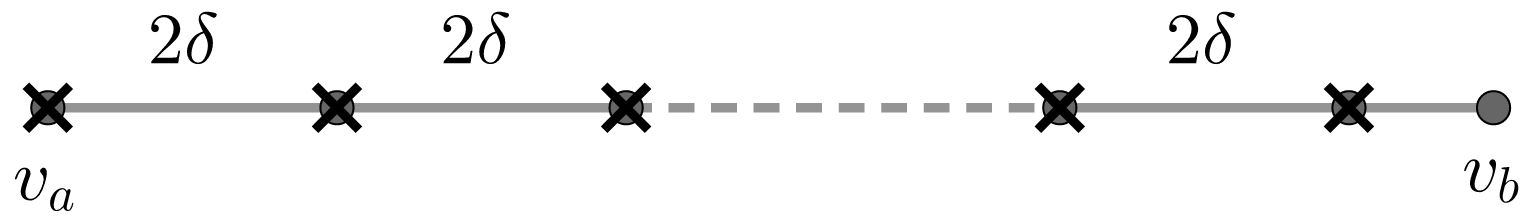}
    	\caption{A facility is located at $v_a$ ($q_e=0$)}
    	\label{fig:brokena}
	\end{subfigure}
	\vspace{0.25cm}
    
 	\begin{subfigure}{0.45\textwidth}
    	\centering
    	\includegraphics[scale=1.1]{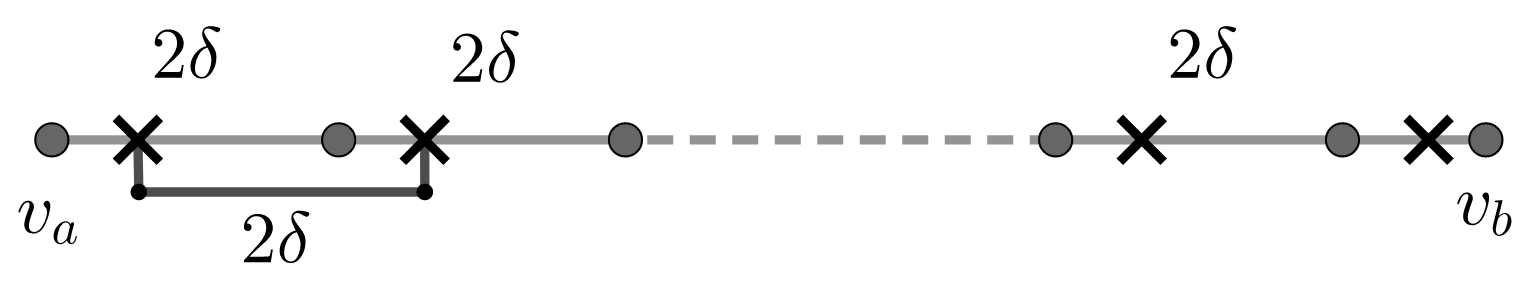}
    	\caption{A facility is located at the tail ($0<q_e\leq \hat{l}_e$)}
    	\label{fig:brokenb}
	\end{subfigure}
	\begin{subfigure}{0.45\textwidth}
    	\centering
    	\includegraphics[scale=1.1]{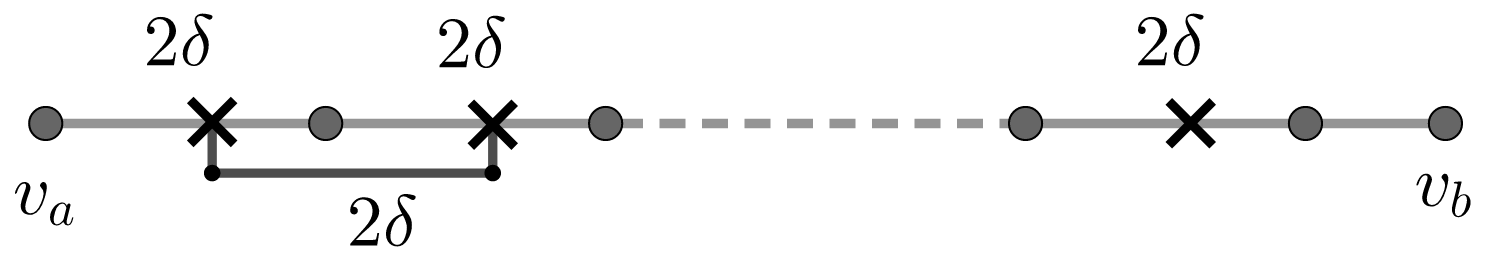}
    	\caption{No facility is located at the tail ($\hat{l}_e <q_e\leq 2\dlt$)}
    	\label{fig:brokenc}
	\end{subfigure}
	\caption{Illustration of Proposition \ref{prop:broken-valid} for a long edge $e$}
	\label{fig:broken}
\end{figure}

The MILP \eqref{milpa} assumes $l_e \le \dlt$ for all $e \in E$, however, in many real-world networks, some edge lengths are greater than the covering radius. To reuse the previous results, one approach is to transform a network with long edges into another network with edge lengths at most $\dlt$. The transformation is by subdividing edges of the original network into smaller pieces, so the optimal cover does not change. This transformation enables us to apply MILP \eqref{milpa} or the MILP in \cite{Hamacher20} on the transformed network. We note that, in \cite{Hamacher20}, there is a recursive transformation. However, the recursive transformation is not constructive, so we cannot employ it in practice.  Another trivial transformation is by subdividing edges as suggested in Section \ref{sec:prelim}.  However, this strategy is a trivial heuristic, and it increases the number of edges and nodes of the transformed network, and thus the number of variables and constraints of the MILP model by a nonlinear factor.

In this section, we present an alternative approach to tackle networks with edge lengths greater than the covering radius. Instead of transforming the network, this approach directly treats long edges in the formulation by using specific sets of constraints and variables. We highlight that the approach is also applicable to ``long paths''. That is, if there is a path in the network whose  intermediate nodes all have degree two, we can represent it by a single edge of length equal to the total length of the path. Indeed, the \problem\ does not change after this transformation.

The main idea of the reduced formulation is to assume a predefined covering of those edges that are long enough. Such covering consists in placing facilities every $2\dlt$ distance units on the long edge. Let us consider $e\in E$ such that $l_e> 2\dlt$. An edge satisfying this condition is called a \emph{long edge}. We denote by $\hat{l}_e:=l_e-2\dlt\lfloor l_e/(2\dlt) \rfloor$ the length of the last piece of $e$ after dividing it into pieces of measure $2\dlt$. We call  $\hat{l}_e$ the \emph{tail} of $e$.  The following proposition guarantees the correctness of the reduced formulation.

\begin{proposition}\label{prop:broken-valid}
Let $N$ be an undirected network, $e=(v_a,v_b)\in E$ be a long edge, and $\cP '$ be a feasible $\dlt$-cover of $N$. Define $\cP$ with $p\in \cP$ for all $p\in \cP ' \setminus e$.  Let $p_e\in \cP '\cap e$ be such that $l(v_a,p_e)=\min_{p'\in \cP '\cap e} l(v_a,p')$, and let $q_e:=l(v_a,p_e)$ (here $q_e$ represents a length, although it will also be a variable of the reduced MILP that we introduce afterward). Note that $q_e\in[0,2\dlt]$ (otherwise, $\cP '$ would not be a $\dlt$-cover). The set $\cP$ can be completed in such a way that it $\dlt$-covers $N$ and $|\cP|\leq | \cP '|$, as follows:
\begin{itemize}
	\item[-] $p_e\in \cP$;
	\item[-] $p\in \cP$ for all $p\in e$, $p>p_e$, such that $l(p_e,p)=2\dlt\cdot k$, for some $k\in \mathbb{N}$;
	\item[-] If $\exists p'\in \cP '\cap e$ such that $d(p_e,p')>2\dlt\lfloor l(p_e,v_b)/(2\dlt) \rfloor$ then $v_b\in \cP $.
\end{itemize}
Moreover, $ \lfloor l_e/(2\dlt) \rfloor\leq |\cP \cap e|\leq \lfloor l_e/(2\dlt) \rfloor +2$. In particular,
\begin{enumerate}
	\item[(i)] If $0\leq q_e\leq \hat{l}_e$, then  $|\cP \cap e|= \lfloor l_e/(2\dlt) \rfloor+1$ if $v_b \notin \cP$,  $|\cP \cap e|= \lfloor l_e/(2\dlt) \rfloor+2$ otherwise.
	\item[(ii)] If $\hat{l}_e <q_e\leq 2\dlt$, then  $|\cP \cap e|= \lfloor l_e/(2\dlt) \rfloor$ if $v_b \notin \cP$,  $|\cP \cap e|= \lfloor l_e/(2\dlt) \rfloor+1$ otherwise.
\end{enumerate}
\end{proposition}
\begin{proof}
It is easy to observe that $|\cP|\leq | \cP '|$. First,  since facilities are placed every 2$\dlt$ distance on $e$, the original covering $\cP '$ cannot contain fewer facilities than $\cP$. Since the rest of the facilities are just taken from $\cP'$, $|\cP|\leq |\cP '|$. On the other hand, $\cP$ has to $\dlt$-cover $N$ for the same reason. That is, the facilities of $\cP$ that were not in $\cP '$ cover at least as much as the ones originally in $\cP'$. The last part of the proposition easily follows from construction, and is illustrated by Figure \ref{fig:broken}.
\end{proof}

 \begin{remark}
The following example shows that the upper bound is tight, that is, there is a case that satisfies $|\mathcal{P} \cap e| =  \lfloor l_e/(2\delta)  \rfloor +2$. Consider a network with a single edge $e$, and let $\delta = 1 $ and $l_e = 3.5$. Therefore, to cover $e$, we need 3 facilities $(|\mathcal{P} \cap e |=3)$. Because $\lfloor l_e/(2\delta)  \rfloor +2 = 3$, the upper bound is tight. Similarly, we can give a tight example for the lower bound, by simply setting $l_e= 2$ in the previous  network.
\end{remark}

In the following, we present our reduced formulation, which is an adaptation of MILP \eqref{milpa}. We treat long edges specifically
, improving the scalability of our approach. Edges $e\in E$ such that $\dlt < l_e \leq 2\dlt$ are subdivided  into two sub edges of length smaller than $\dlt$. Therefore, we assume that, for every $e\in E$, either $l_e\leq \dlt$ or $l_e> 2\dlt$. In the former case, all the constraints and variables of the model remain unchanged. In the latter, we introduce new variables and constraints to the model, while dropping some of the constraints originally in \eqref{milpa}. The objective function also needs adaptation. We introduce all these modifications next.

Let $e=(v_a,v_b)\in E$ be a long edge and, for any feasible solution, let $q_e$ be as in Proposition \ref{prop:broken-valid}. That is, when $e$ is a long edge, we use the former variable $q_e$ of MILP \eqref{milpa} to represent the position of the left-most facility on $e$ with respect to $v_a$. The placement variables $y_{v_a}$, $y_{v_b}$, and $y_e$ will be used as well, with slightly different meanings to those in \eqref{milpa}, as we will explain later on. We introduce an indicator variable $u_{e}\in\{0,1\}$ to distinguish between two possible ranges in the domain of $q_e$.  If $0 \le q_{e} \le \hat l_e$, then $u_{e} =0$; otherwise $\hat l_{e} \le q_{e} \le 2\dlt$ and $u_{e} =1$.  This can be modeled with the following constraints:
\begin{equation}
\label{eq.cvu}
\begin{split}
	q_{e} &\le \hat l_e (1-u_{e})   +2\dlt u_{e},\\
	q_{e} &\ge \hat l_e u_{e}. \\
\end{split}
\end{equation}
From Proposition \ref{prop:broken-valid},  there is a  transition in the number of facilities on $e$ when $u_{e}$ changes from $0$ to $1$. Let us denote by $L\subseteq E$ the set of long edges of the network. The objective function of the reduced MILP reads as follows
\begin{equation}\label{obj-new}
	\sum_{f\in\cF\setminus L} y_f + \sum_{e\in L} \Big(\left\lceil \frac{l_e}{2\dlt}\right\rceil -u_e\Big).
\end{equation}
Note that the coefficients on the last term in the objective already account for the facilities installed at $v_a$ for each $e=(v_a,v_b)\in L$, while they do not do so for $v_b$. This will condition the values of the placement variables in an optimal solution, namely, $y_{v_a}=0$ for all $e=(v_a,v_b)\in L$. The facilities installed at these nodes will be tracked by the variables $q_e$, namely, if $q_e=0$ then a facility would be installed at $v_a$. To complete the modeling of the \problem, we need to ensure that the covering of the long edge fits into the covering of the rest of the network. Namely, some parts of the network might be covered by facilities placed on $e$, and part of $e$ (namely its tail or the portion between $v_a$ and $p_e$ of Proposition \ref{prop:broken-valid}) could be covered by facilities placed outside $e$.

We focus first on the case $0 \le q_{e} \le  \hat l_{e}$. We need to ensure that both the segment $(v_a,p_e)$ and the tail of $e$ are covered.
For the tail, we know that there is a facility at a distance $\hat l_e-q_{e}$ from $v_b$. This facility covers a length $\dlt$ of the remaining fragment on $e$ on its right-hand side, which has a length equal to $ \hat l_{e}-q_{e}$.  The rest of such fragment should be covered, which can be imposed by the following constraint:
  \begin{equation}
  \label{eq.cvu0}
  	r_{v_{b}} \ge \hat l_e-q_{e}-\dlt \iff	r_{v_{b}}+q_{e} + \dlt \ge \hat l_e \qquad \forall e \in L \st u_{e} =0.
  \end{equation}
To ensure the covering of the segment $(v_a,p_e)$, we have:
	\begin{equation}
 	r_{v_{a}}+\dlt \ge q_{e} \qquad  \forall e \in L \st u_{e} =0. \nonumber
  \end{equation}
 
 
Let us consider now the case $\hat l_{e} < q_{e} \le 2\dlt$.  There is a facility installed at a distance of $2\dlt - q_{e} + \hat l_e$ from $v_b$. Then, to ensure that the tail of $e$ is covered, we need to cover the fragment between this facility and $v_b$. Since the facility already covers a length $\dlt$ on this fragment, the following constraint enforces the covering of the tail:
  \begin{equation}
  \label{eq.cvu1}
  	r_{v_{b}}\ge 2\dlt  - q_{e}+ \hat l_{e}  -\dlt   \iff  r_{v_{b}}+q_{e} - \dlt \ge\hat l_{e}   \qquad\forall e \in L \st u_{e} =1.
  \end{equation}
To ensure the covering of the segment $(v_a,p_e)$, we have the same equation as before:
	\begin{equation}
 	r_{v_{a}}+\dlt \ge q_{e} \qquad  \forall e \in L \st u_{e} =1. \nonumber
  \end{equation}

In summary, the following constraints ensure that the edge $e$ is fully covered:
\begin{eqnarray}
	r_{v_{a}}+\dlt \ge q_{e} & \forall e\in L,\label{eq.cvhead}\\
	r_{v_{b}} + q_{e} - (2u_{e}-1)\dlt \ge \hat l_e & \forall e\in L.\label{eq.cvtail}
\end{eqnarray}
Constraint \eqref{eq.cvtail} gathers \eqref{eq.cvu0} and \eqref{eq.cvu1} in a single constraint. The reduced MILP model is as follows (we avoid extended writing of the model for the sake of conciseness):
\begin{eqnarray}
	\min & \eqref{obj-new} \nonumber\\
	\st   & \eqref{milpa.completecover1}, \eqref{milpa.completecover2}, \eqref{milpa.covere} & e\notin L \nonumber\\
	& \eqref{milpa.optimale}, \eqref{milpa.coord},\eqref{milpa.indicatorei} & e'\notin L \nonumber\\
	&  \eqref{milpa.enforcex1}, \eqref{milpa.enforcex2}, \eqref{milpa.lcoversos},\eqref{milpa.indicatorv},\eqref{milpa.wbdl},\eqref{milpa.coverdist-node},\eqref{milpa.coverdist-edge}, \eqref{milpa.varyw},\eqref{milpa.varindicator},\eqref{milpa.varcontinuous}\nonumber\\
	& \eqref{eq.cvu}, \eqref{eq.cvhead}, \eqref{eq.cvtail} \nonumber\\
	& y_{e}=1 & e \in L \label{fix-y}\\
	& w_{e}=0 & e \in L, \label{fix-w}
\end{eqnarray}
 where, if $e\in L$, the term $\tau_{vei}(q_{e})$ in \eqref{milpa.coverdist-edge} is replaced by $d(v,v_{i})+ \mathbf{1}_{i=a}q_{e} +\mathbf{1}_{i=b}(2\dlt u_{e} + \hat{l}_e  - q_{e})$.
 We enforce \eqref{fix-y} because $e$ always contains a facility if $e\in L$. On the other hand, we need to include constraints \eqref{fix-w} to guarantee that the variables $r_{v_a}$ and $r_{v_b}$ can take positive values. Indeed, if $w_e=1$ for $e\in L$, it may happen that $x_{v_a}=1$ or $x_{v_b}=1$, which will will imply, respectively,  $r_{v_a}=0$ or $r_{v_b}=0$  due to \eqref{milpa.wbdl}.
We compute the complete and partial cover sets in the same way as  for the original MILP model. Note that no edge or node can completely cover $e$ if $e\in L$.

The following theorem is on the scalability of the reduced MILP above.
\begin{theorem}
Given a network 
$N=(V,E,l)$, the maximum number of variables and constraints of the reduced MILP model only depends on $V$ and $E$.
\end{theorem}
\begin{proof}
The number of constraints and variables of the reduced model does not grow with the edge lengths, except for a constant factor of 2 for those edges $e\in E$ such that $\dlt < l_e \leq 2\dlt$.
\end{proof}

\section{Computational results}\label{sec:compu}

In this section, we present the computational experiments testing the existing and proposed formulations and strengthening techniques for  \problem \, and its discrete variant (facilities on nodes).

\subsection{Experiment Setup}
We describe the setup of the experiments including the benchmarks, development environment, implementation of algorithms and solution statistics. The  computational results and source code are publicly released on our project website: \href{https://github.com/lidingxu/cflg/}{https://github.com/lidingxu/cflg/}, where we provide a bash file to reproduce the experiments in Linux systems. Those benchmarks that we generated for this study, or that were publicly available already, are also available at the repository. 

\textbf{Benchmarks.}
We use three different benchmarking sets: two come from the literature, and the other has been generated synthetically. For every instance, we set the coverage radius $\dlt$ equal to the  average of the edge lengths. We describe these benchmarks next.

\texttt{Kgroup.} It consists of 23 prize-collecting Steiner tree problem instances from \cite{Ljubic2006}, and the benchmark includes the graphs and edge lengths of these instances. These random geometric instances are designed to have a local structure somewhat similar to street maps. Nodes correspond to random points in the unit square. The number of nodes ranges from 22 to 241. There
is an edge between two nodes if their distance is no more than a prescribed threshold which depends on the number of nodes, and the length of an edge is the Euclidean distance between the two points. It is divided into two sets, \texttt{Kgroup\_A} and  \texttt{Kgroup\_B}. The first one consists of 12 small instances with up to 45 nodes, and the second one consists of 11 large instances with up to 241 nodes.  

\texttt{City.} It consists of real data of 9  street networks for some German cities, and it was first used in \cite{Kalsics}.  The number of nodes ranges from 132 to 771. The length of each edge is the length of the underlying street segment.

\texttt{Random.} It consists of 24  random network instances generated via Erdős-Rényi binomial method with the package ``Networkx" (see \cite{SciPyProceedings_11}).  A network is constructed by connecting nodes randomly. Each edge is included with a predefined uniform probability $p$. The number of nodes, $n$, is in $\{10, 15, 20, 25, 30, 40\}$. For each $n$, we generate random graphs with different adjacency probabilities, namely  $p \in \{0.1, 0.2, 0.3, 0.4\}$. Furthermore, we split these instances into two benchmarks: \texttt{Random\_A} and \texttt{Random\_B}.  \texttt{Random\_A} contains instances with $n \in \{10,15,20\}$.   \texttt{Random\_B} contains instances with $n \in \{25,30.40\}$. \\

\textbf{Coverage radii.}
For each network, we define two sets of coverage radii: ``Small'' equal to [Average Edge
Length], and ``Large'' equal to ×2 [Average Edge Length], respectively.
\\


\textbf{Problem preprocessing.}
Networks of instances are modified in a problem preprocessing step to be amenable to MILP models.

 Given an original network of each instance, in the first preprocessing step, we delete any degree-two node and concatenate its adjacent edges to a new edge, as long as the deletion does not yield a self-loop. Such a node can be treated as an interior point of the new edge. We refer to the preprocessed network without any such degree-two node as the degree-two-free network.

Even after the first preprocessing step,
 the degree-two-free network may not correspond to the actual problem network to solve, since we may subdivide the degree-two-free network for the non-reduced model to guarantee that $\dlt > \max_{e \in E} |e|$. We refer to the preprocessed network after the second preprocessing step as the subdivided network, which is degree-two-free and satisfies $\dlt > \max_{e \in E} |e|$.
 Therefore, the size (number of nodes and edges) of a subdivided network depends on $\dlt$.\\

\textbf{Development environment.}
The experiments are conducted on a computer with Intel Core i7-6700K  CPU @ 4.00GHZ and 16GB main memory. JuMP \cite{DunningHuchetteLubin2017} is a modeling language for mathematical optimization embedded in Julia. We use JuMP  to implement our models and interact with MILP solvers. Specifically, we use ILOG CPLEX 20.1 to solve our models. Alternatively, the implementation allows users to switch easily to other solvers (e.g. Gurobi and GLPK).

CPLEX's parameters are set as their defaults, except that we disable its parallelism and set  the MIP absolute gap to 1 (due to the integral objective). The experiments are partitioned into jobs. Every job calls CPLEX to solve an instance, and this job is handled by one process of the multi-core CPU. To safeguard
against a potential mutual slowdown of parallel processes, we run only one job per core at a
time, and we use at most three processes in parallelism.  The time limit of each job is set to 1800 CPU seconds.\\

\textbf{Model implementation.}
We implement six models based on different combinations of formulations and settings. The first five models address \problem, while the last model solves its discrete restriction, i.e. the variant in which facilities must be placed at nodes. These models are as follows.

\texttt{EF.} This model implements the model from \cite{Hamacher20} for \problem. This formulation only uses edges to model facility locations, and the authors do not consider the complete and partial cover sets to delimit the size of the model. This model assumes $\dlt > \max_{e \in E}|e|$, and it reads the subdivided graph.

\texttt{F0.} This model implements a  basic formulation that is a simplification of the model \eqref{milpa}. It does not use the complete and partial cover information nor any of the strengthening techniques in \Cref{sec:strenght}. Hence, it does not call the network processing algorithm \ncalgo. This model assumes $\dlt > \max_{e \in E}|e|$, and it reads the subdivided graph. More precisely, the constraints \eqref{milpa.completecover1}-\eqref{milpa.enforcex2} related to complete covers are removed, the complete cover variables $w$ are fixed to $0$; for each $v \in V$, the partial cover sets $\cE_{\mathrm{p}}(v)$, $\cEI_{\mathrm{p}}(v)$ are solely set, respectively, as $E$ and $E \times \{a,b\}$, and consequently, $M_{v} = \dlt$, $M_{vv'} = r(N)$ for $v' \in \cE_{\mathrm{p}}(v)$, $M_{ve'i'} = r(N) + |e'|$ for $(e',i') \in \cEI_{\mathrm{p}}(v)$ are trivial valid bound constants, where $r(N):= \max_{v,v' \in N} d(v,v')$ is the radius of the problem network $N$.

\texttt{F.}  This model implements the complete formulation \eqref{milpa} for \problem, it does use the complete and partial cover information, and hence it  calls the network processing algorithm \ncalgo.  It does not use the strengthening techniques in \Cref{sec:strenght}. This model assumes $\dlt > \max_{e \in E}|e|$, and it reads the subdivided network as well.  For each $v \in V$, due to the delimited partial cover set,  $M_{v} = \dlt$, $M_{vv'} = \dlt$ for $v' \in \cE_{\mathrm{p}}(v)$, $M_{ve'i'} = \dlt + |e'|$ for $(e',i') \in \cEI_{\mathrm{p}}(v)$ are valid bound constants.

\texttt{SF.} This model strengthens  $\texttt{F}$ by using the techniques described in \Cref{sec:strenght}. 
More precisely, the big-M constants are reduced as \Cref{sec:strenght.bd}; the "Leafs" inequalities are used to fix variables; and the "Neighborhood" inequalities are implemented as model constraints which replace \eqref{milpa.optimale}.

\texttt{RF.} This model implements the reduced formulation from \Cref{sec:reduced-milp}. It only requires $\dlt < 2\max_{e \in E}|e|$. Given a degree-two-free network, it models the long edge specifically as the description \Cref{sec:reduced-milp}, and it subdivides the edges with lengths greater than $\dlt$ and smaller than $2\dlt$ into two sub-edges.

\texttt{SFD.}  Any solution of the discrete restriction of \problem---where facilities can only be placed at nodes--- is a feasible solution of \problem. We name this discrete restriction by the \dfproblem.  This model solves the \dfproblem, which solely sets $y_e = 0$ for all $e \in E$ in \texttt{SF} model.

The above models are summarized in \Cref{tab.modsum}. Both \texttt{EF} and \texttt{F0} consider that any two points in the network can possibly cover each other, and do not utilize the complete and partial cover information. They have been already compared in \Cref{sec:comp}, and hence \texttt{F0} should have fewer variables and constraints than \texttt{EF}. We are interested in the dual gaps obtained after the models are solved within the time limit for these models. \\

 \begin{table}[]
 \resizebox{0.99\columnwidth}{!}{
\begin{tabular}{|l|lllllll|}
\hline
  \multicolumn{1}{|l|}{Model}   & \multicolumn{1}{l|}{Problem} & \multicolumn{1}{l|}{Delimitation} & \multicolumn{1}{l|}{Strengthening} & \multicolumn{1}{l|}{Long edge} & \multicolumn{1}{l|}{Size} & \multicolumn{1}{l|}{Input network} & \multicolumn{1}{l|}{Comment} \\ \hline
\texttt{EF}  & \problem                                  	& No                                	& No                              	& No                         	& Very large                  	& Subdivided network             	&            	From \cite{Hamacher20}          	\\ \cline{1-8}
\texttt{F0}  &  \problem                                    	& No                                	& No                              	& No                         	& Large                       	& Subdivided network             	& The simple model          	\\ \cline{1-8}
\texttt{F}   &  \problem                                    	& Yes                               	& No                              	& No                         	& Medium                      	& Subdivided network             	& The complete model       	\\ \cline{1-8}
\texttt{SF}  &  \problem                                    	& Yes                               	& Yes                             	& No                         	& Medium                      	& Subdivided network             	& The strengthened model   	\\ \cline{1-8}
\texttt{RF}  &  \problem                                    	& Yes                               	& Yes                             	& Yes                        	& Small                       	& Degree-two-free network        	& The reduced model        	\\ \cline{1-8}
\texttt{SFD} &  \dfproblem                                    	& Yes                               	& Yes                             	& No                         	& Very Small                  	& Subdivided network             	& The discrete model       	\\ \cline{1-8}
\end{tabular}}
\caption{Model summary}
\label{tab.modsum}
\end{table}

\textbf{Performance metrics and statistical tests.}
We describe the performance metrics and the ways to compute  their statistics. These statistics will be used to evaluate the model performance.

Let $\underline{v}$ be a dual lower bound and $\overline{v}$ be a primal upper bound obtained after solving some of the models described above,  the relative dual gap   is defined as:
\begin{equation*}
   \sigma := \frac{\overline{v} - \underline{v}}{\overline{v}}.
\end{equation*}
A smaller relative dual gap indicates better primal and dual behavior of the model.
 
Let $n_{sd}$ be the number of nodes of the subdivided network of that instance, note that  a trivial primal solution is the set of the nodes of the subdivided network (for which edge length is at most $\dlt$). Therefore, to normalize the primal solution value, we define the relative primal bound  $$v_r := \frac{\overline{v}}{n_{sd}} .$$ If $v_r < 1$, then the model finds a solution better than the trivial one.

In order to evaluate model performance, we compute shifted geometric means (SGMs) of performance metrics, which provides a
measure for relative differences. This avoids statistics from being dominated by outliers with
large absolute values as is the case for the arithmetic mean. The SGM also avoids an over-representation of results with small absolute values.  The SGM of values $v_1,...,v_M \geq 0$ with shift $s \geq 0$ is defined
as
\begin{equation*}
  \left(\prod_{i=1}^M (v_i + s)\right)^{1/M} - s.
\end{equation*}

We say an instance is affected by a model, if solving this model finds a feasible solution; the instance is solved by this model, if solving this model finds an optimal solution. If an instance is unaffected, usually the model is too large to be read into the MILP solver.

We record the following performance metrics of each instance for each model, and compute the benchmark-wise SGMs:

\begin{enumerate}
	\item $t$: the total running time in CPU seconds, with a shifted value set to 1 second;
    
   \item $\sigma$: the relative dual gap, with a shifted value set to $0.01$;
  	\item $v_r$: the relative primal bound, with a shifted value set to $0.01$.
	 
\end{enumerate}
For an unaffected instance, we set by default $t= 1800$, $\sigma = 1$ and $v_r = 1$. Note that the time does not include the preprocessing time, since we find that the preprocessing is usually at most 0.5 seconds. \\

We will discuss the computational results, which are divided into two parts. In the first part, we compare the five models \texttt{EF}, \texttt{F0}, \texttt{F}, \texttt{SF}, and \texttt{RF}. We evaluate the performance metrics of these models. The second part compares \texttt{RF}  and \texttt{SFD},  quantifying the facilities that are saved by allowing continuous location. In the following, we will analyze the aggregated results. We refer to  \Cref{table:root:detailed}  in
the appendix for the detailed instance-wise results, where decimal representations of gaps and bounds are converted to percentages.

\subsection{Comparative Analysis of Continuous Models}

We compare five continuous models for the \problem, namely  \texttt{EF}, \texttt{F0}, \texttt{F}, \texttt{SF}, and \texttt{RF}. For each  benchmark, radius and model, we record a triple of integers S/A/T: S denotes the number of solved instances, A denotes the number of affected instances, and T denotes the number of total instances in the benchmark. Moreover, we also report the  average SGMs of the dual gaps, solving times and relative primal bounds among all instances in the benchmark. \Cref{tab.cm} summarizes these results.
 
 First, we notice that \texttt{EF} cannot affect any instance in any benchmark; \texttt{RF}, \texttt{SF} and \texttt{RF} can affect all instances, i.e., solutions are provided by these models; \texttt{RF} is the model that solves the most number of instances (11), and \texttt{SF} is the second best one (10).
 
 Secondly, we compare \texttt{EF} and \texttt{F0}. \texttt{F0} is obviously superior to \texttt{EF}. With \texttt{F0}, 39 among 56 instances of small radius (resp. 42 among 56 instances of large radius) can be read by the CPLEX solver, while the instances modeled by \texttt{EF} are too large to read. Therefore, better solutions than trivial solutions are found by \texttt{F0}: on average, for instance of small radius (resp. large radius), \texttt{F0} finds solutions that use $25.2\%$ (resp. $74.1\%$) fewer facilities than the trivial solution.

Then, we compare \texttt{F0}, \texttt{F}  and \texttt{SF}. With the delimitation of complete and partial covering sets, \texttt{F} and \texttt{SF} can affect all instances (especially those in \texttt{Kgroup\_B}, of which \texttt{F0} could just read one).
With the strengthening technique, \texttt{SF} has only marginal improvement in the relative primal bound, and solving time, while \texttt{F} is even slightly better than \texttt{SF} in the dual gap. We observe, in our experiment, that adding valid inequalities might slow down the internal solving process of CPLEX.

Finally, we compare  \texttt{SF} and \texttt{RF}. \texttt{RF} outperforms \texttt{SF} in all performance metrics. Moreover,  \texttt{RF} is the best one among those models affecting all instances. Indeed, for many instances, their degree-two-free networks may contain long edges, and \texttt{RF} avoids introducing too many variables and constraints for modeling their coverage.

\texttt{Kgroup\_B} is the hardest benchmark. The best model \texttt{RF} still has  an average dual gap of $59.1\%$ and $154.2\%$ relative primal bound for instances of small radius, and this means that $\texttt{RF}$ cannot produce better solutions than the trivial one.

We find that for all the models (except for \texttt{EF}), the average dual gaps and solving times of instances of large radius are smaller than those of instances of small radius. This shows that the large radius has a positive effect on the model performance, and an instance of a small radius may be more difficult than the same instance of a large radius. This is because, with a larger radius, the network after processing is smaller (see \Cref{table:root:detailed}).

\begin{figure}[!h]
\centering{
\includegraphics[width = 0.999\columnwidth]{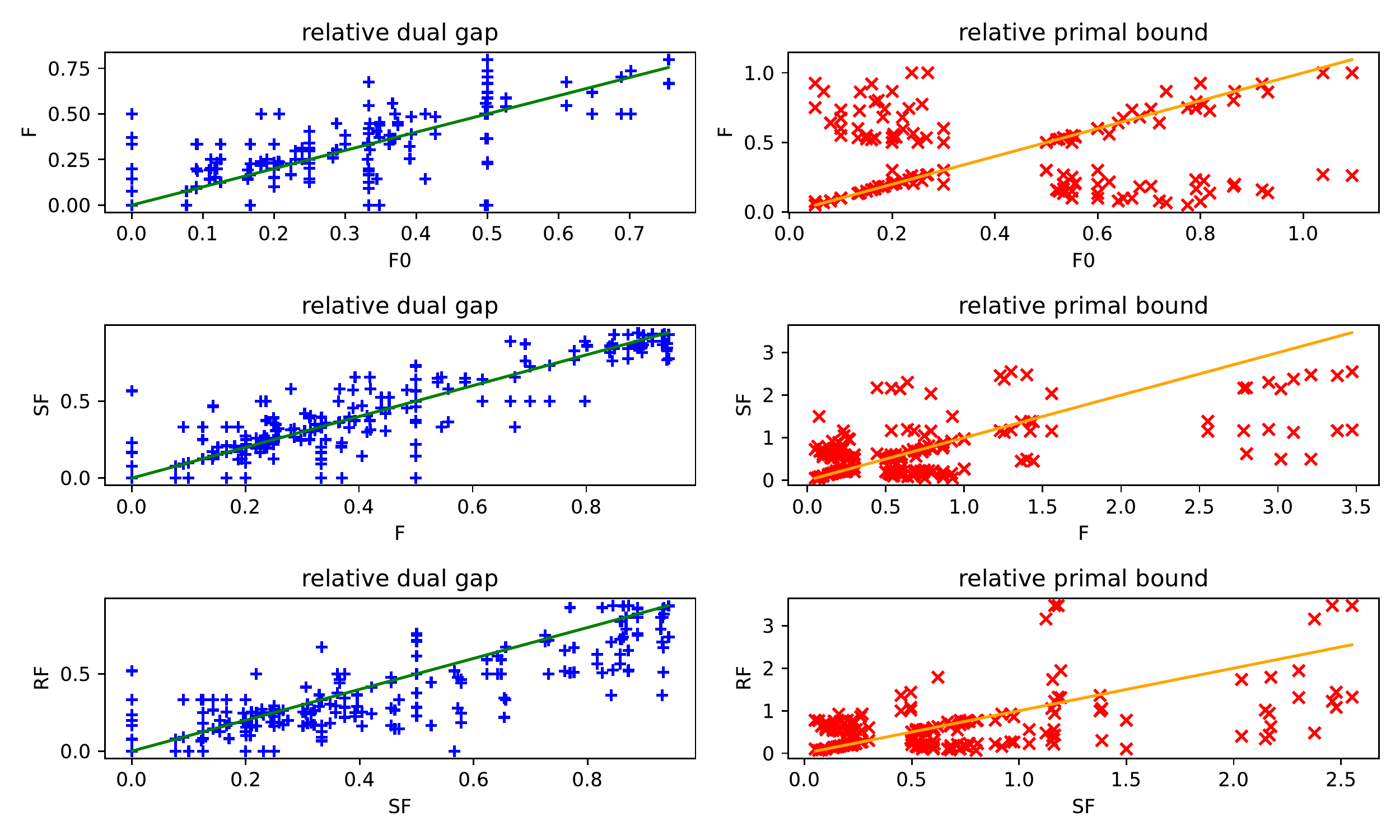}
}
\caption{Scatter plots of the relative dual gaps and the relative primal bounds between different settings}
\label{fig.scatter5}
\end{figure}

In \Cref{fig.scatter5}, we show scatter plots of the relative dual gaps and the relative primal bounds  of affected instances between different settings. For every plot, there is a line in which the points have equal (X,Y)-values. If points fall below the line, then the Y-axis model performs better for the corresponding instances. Note that when comparing  \texttt{F0} and  \texttt{F}, the plots do not consider the unaffected instances of \texttt{F0} which are affected or solved by \texttt{F}. Moreover, \texttt{F0} even closes more duality gaps than \texttt{F}, but \texttt{F} can find better primal solutions. These plots give an overview of all affected instances and support the above analysis.

To summarize, we have shown that the two proposed techniques---  that to delimit the coverage  areas from a given point in Section \ref{sec:prelim}, and that to cover long edges in Section \ref{sec:reduced-milp},--- can reduce the model size drastically. Among the five models tested, $\texttt{RF}$ features the best overall performance, which is achieved by directly modeling covers on long edges. On the other hand, delimiting the covering sets to the potential, complete and partial covers also reduce the model size, which allows $\texttt{F}$ to read all the tested instances.

\begin{table}[]
\centering

\begin{tabular}{|l|l|*{4}{r}|*{4}{r}}
\hline
\multirow{2}{*}{Benchmark} & \multirow{2}{*}{Radius} & \multicolumn{4}{l|}{\texttt{EF}}                                                                                 	& \multicolumn{4}{l|}{\texttt{F0}}                                                                                                    	\\ \cline{3-10}
                       	& & \multicolumn{1}{l|}{time} & \multicolumn{1}{l|}{$\sigma(\%)$} & \multicolumn{1}{l|}{$v_r(\%)$} & \multicolumn{1}{l|}{S/A/T} & \multicolumn{1}{l|}{time} & \multicolumn{1}{l|}{$\sigma(\%)$} & \multicolumn{1}{l|}{$v_r(\%)$} & \multicolumn{1}{l|}{S/A/T}  \\ \hline
	\multirow{2}{*}{\texttt{city}}&Small & 1800.0 & 100.0\% & 100.0\% & 0/0/9 & 1801.7 & 56.8\% & 83.3\% & \multicolumn{1}{l|}{0/3/9} \\
&Large & 1800.0 & 100.0\% & 100.0\% & 0/0/9 & 1800.9 & 42.3\% & 36.2\% & \multicolumn{1}{l|}{0/6/9} \\ \cline{1-1}
\multirow{2}{*}{\texttt{Kgroup\_A}}&Small & 1800.0 & 100.0\% & 100.0\% & 0/0/11 & 1802.6 & 25.1\% & 85.0\% & \multicolumn{1}{l|}{0/11/11} \\
&Large & 1800.0 & 100.0\% & 100.0\% & 0/0/11 & 139.2 & 14.7\% & 19.2\% & \multicolumn{1}{l|}{7/11/11} \\ \cline{1-1}
\multirow{2}{*}{\texttt{Kgroup\_B}}&Small & 1800.0 & 100.0\% & 100.0\% & 0/0/12 & 1800.4 & 92.6\% & 98.8\% & \multicolumn{1}{l|}{0/1/12} \\
&Large & 1800.0 & 100.0\% & 100.0\% & 0/0/12 & 1800.1 & 93.2\% & 86.6\% & \multicolumn{1}{l|}{0/1/12} \\ \cline{1-1}
\multirow{2}{*}{\texttt{random\_A}}&Small & 1800.0 & 100.0\% & 100.0\% & 0/0/12 & 16.8 & 15.9\% & 54.8\% & \multicolumn{1}{l|}{9/12/12} \\
&Large & 1800.0 & 100.0\% & 100.0\% & 0/0/12 & 0.2 & 25.5\% & 19.5\% & \multicolumn{1}{l|}{12/12/12} \\ \cline{1-1}
\multirow{2}{*}{\texttt{random\_B}}&Small & 1800.0 & 100.0\% & 100.0\% & 0/0/12 & 1317.6 & 36.4\% & 63.3\% & \multicolumn{1}{l|}{1/12/12} \\
&Large & 1800.0 & 100.0\% & 100.0\% & 0/0/12 & 154.4 & 26.0\% & 10.0\% & \multicolumn{1}{l|}{11/12/12} \\ \cline{1-10}
\multirow{2}{*}{\texttt{all}}&Small & 1800.0 & 100.0\% & 100.0\% & 0/0/56 & 625.8 & 37.4\% & 74.8\% & \multicolumn{1}{l|}{10/39/56} \\
&Large & 1800.0 & 100.0\% & 100.0\% & 0/0/56 & 132.5 & 33.1\% & 25.9\% & \multicolumn{1}{l|}{30/42/56} \\
 \cline{1-10}
\multicolumn{10}{c}{}\\
\cline{1-10}
\multirow{2}{*}{Benchmark} & \multirow{2}{*}{Radius} & \multicolumn{4}{l|}{\texttt{F}}                                                                                 	& \multicolumn{4}{l|}{\texttt{SF}}                                                                                                    	\\ \cline{3-10}
                       	& & \multicolumn{1}{l|}{time} & \multicolumn{1}{l|}{$\sigma(\%)$} & \multicolumn{1}{l|}{$v_r(\%)$} & \multicolumn{1}{l|}{S/A/T} & \multicolumn{1}{l|}{time} & \multicolumn{1}{l|}{$\sigma(\%)$} & \multicolumn{1}{l|}{$v_r(\%)$} & \multicolumn{1}{l|}{S/A/T}  \\ \hline
\multirow{2}{*}{\texttt{city}}&Small & 1802.9 & 29.5\% & 62.2\% & 0/9/9 & 1801.3 & 30.1\% & 66.9\% & \multicolumn{1}{l|}{0/9/9} \\
&Large & 1801.2 & 28.4\% & 21.7\% & 0/9/9 & 1800.9 & 29.1\% & 21.7\% & \multicolumn{1}{l|}{0/9/9} \\  \cline{1-1}
\multirow{2}{*}{\texttt{Kgroup\_A}}&Small & 1803.0 & 33.1\% & 82.2\% & 0/11/11 & 1801.3 & 32.0\% & 80.6\% & \multicolumn{1}{l|}{0/11/11} \\
&Large & 238.0 & 18.9\% & 19.1\% & 8/11/11 & 300.8 & 19.0\% & 19.1\% & \multicolumn{1}{l|}{8/11/11} \\  \cline{1-1}
\multirow{2}{*}{\texttt{Kgroup\_B}}&Small & 1800.6 & 80.8\% & 240.5\% & 0/12/12 & 1801.4 & 79.7\% & 191.9\% & \multicolumn{1}{l|}{0/12/12} \\
&Large & 1800.4 & 85.1\% & 80.5\% & 0/12/12 & 1800.7 & 85.9\% & 77.3\% & \multicolumn{1}{l|}{0/12/12} \\  \cline{1-1}
\multirow{2}{*}{\texttt{random\_A}}&Small & 20.2 & 16.5\% & 54.3\% & 9/12/12 & 16.1 & 17.1\% & 54.9\% & \multicolumn{1}{l|}{9/12/12} \\
&Large & 0.3 & 25.5\% & 19.5\% & 12/12/12 & 0.2 & 10.4\% & 17.9\% & \multicolumn{1}{l|}{12/12/12} \\  \cline{1-1}
\multirow{2}{*}{\texttt{random\_B}}&Small & 1574.2 & 38.8\% & 64.9\% & 1/12/12 & 1501.2 & 40.0\% & 67.5\% & \multicolumn{1}{l|}{1/12/12} \\
&Large & 220.5 & 19.9\% & 10.3\% & 9/12/12 & 175.7 & 18.8\% & 10.0\% & \multicolumn{1}{l|}{11/12/12} \\  \cline{1-10}
\multirow{2}{*}{\texttt{all}}&Small & 675.0 & 35.2\% & 86.2\% & 10/56/56 & 637.6 & 35.5\% & 83.6\% & \multicolumn{1}{l|}{10/56/56} \\
&Large & 163.0 & 30.2\% & 23.6\% & 29/56/56 & 160.9 & 24.9\% & 22.8\% & \multicolumn{1}{l|}{31/56/56} \\

\cline{1-10}
\multicolumn{10}{c}{}\\
\cline{1-6}
\multirow{2}{*}{Benchmark} & \multirow{2}{*}{Radius}  & \multicolumn{4}{l|}{\texttt{RF}}                                                                                 	& \multicolumn{4}{l}{}                                                                                                    	\\ \cline{3-6}
                  	&  	& \multicolumn{1}{l|}{time} & \multicolumn{1}{l|}{$\sigma(\%)$} & \multicolumn{1}{l|}{$v_r(\%)$} & \multicolumn{1}{l|}{S/A/T} & \multicolumn{1}{l}{} & \multicolumn{1}{l}{} & \multicolumn{1}{l}{} & \multicolumn{1}{l}{}  \\ \cline{1-6}
\multirow{2}{*}{\texttt{city}}&Small & 1804.4 & 16.2\% & 54.1\% & \multicolumn{1}{l|}{0/9/9} \\
&Large & 1801.5 & 25.8\% & 21.3\% & \multicolumn{1}{l|}{0/9/9} \\  \cline{1-1}
\multirow{2}{*}{\texttt{Kgroup\_A}}&Small & 1622.6 & 21.5\% & 77.5\% & \multicolumn{1}{l|}{1/11/11} \\
&Large & 158.9 & 19.2\% & 19.3\% & \multicolumn{1}{l|}{8/11/11} \\  \cline{1-1}
\multirow{2}{*}{\texttt{Kgroup\_B}}&Small & 1800.9 & 59.1\% & 154.2\% & \multicolumn{1}{l|}{0/12/12} \\
&Large & 1800.6 & 75.5\% & 63.3\% & \multicolumn{1}{l|}{0/12/12} \\  \cline{1-1}
\multirow{2}{*}{\texttt{random\_A}}&Small & 15.9 & 8.1\% & 54.3\% & \multicolumn{1}{l|}{9/12/12} \\
&Large & 0.3 & 26.6\% & 19.8\% & \multicolumn{1}{l|}{12/12/12} \\  \cline{1-1}
\multirow{2}{*}{\texttt{random\_B}}&Small & 1304.3 & 38.5\% & 63.8\% & \multicolumn{1}{l|}{1/12/12} \\
&Large & 190.2 & 19.8\% & 11.2\% & \multicolumn{1}{l|}{9/12/12} \\  \cline{1-6}
\multirow{2}{*}{\texttt{all}}&Small & 604.9 & 23.7\% & 75.4\% & \multicolumn{1}{l|}{11/56/56} \\
&Large & 146.6 & 29.2\% & 22.8\% & \multicolumn{1}{l|}{29/56/56} \\
\cline{1-6}
\end{tabular}

\caption{Results for continuous models} \label{tab.cm}
\end{table}

\subsection{Comparative Analysis of Continuous and Discrete Models}

 In  \problem, the facilities are located either at nodes or edges, while in the discrete variant considered in this section facilities can only be located at nodes.
Our objective is to evaluate the number of facilities that can be saved by allowing continuous location. Since the discrete model studied here, \texttt{SFD}, is a discrete restriction of \problem,  every optimal solution is a feasible solution of \problem.
 We solve \texttt{SFD} for the \dfproblem \,  and compare the results with the best model for \problem, \texttt{RF}.

In addition to the previous performance statistics, we also record for each instance, a new relative primal bound for the continuous model defined as:
$$v'_r := \frac{\overline{v}}{\overline{v_d}},$$ where $\overline{v_d}$ is the best solution found by \texttt{SFD}.  If $v'_r < 1$, then the continuous model (in this case, \texttt{RF}) finds a solution better than the one found by the discrete model.

 \begin{figure}[htbp]
\centering{
\includegraphics[width = 0.999\columnwidth]{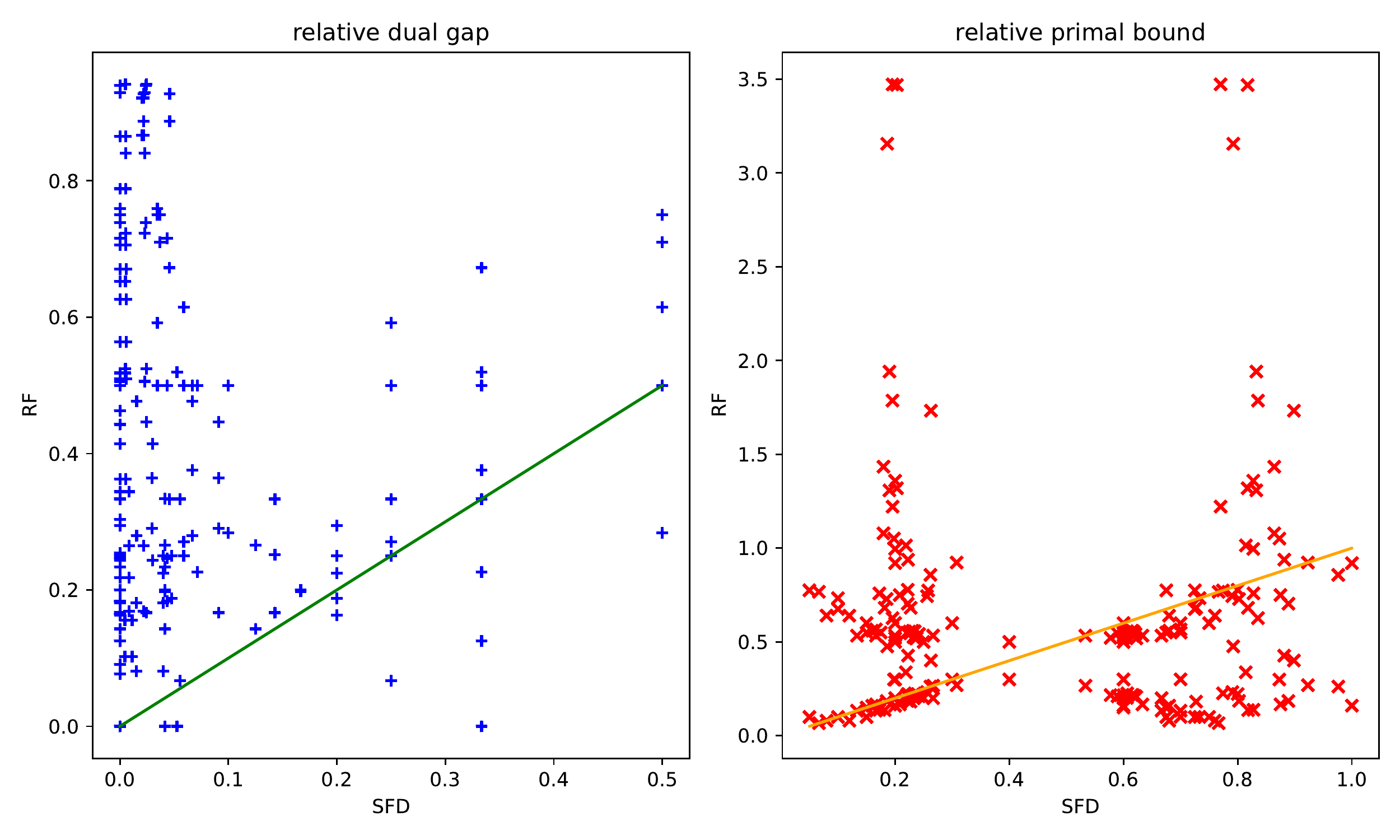}
}
\caption{Scatter plots of the relative dual gaps and the relative primal bounds between \textbf{SFD} and \textbf{RF}}
\label{fig.scatter2}
\end{figure}

\begin{table}[htbp]
\begin{center}
\begin{tabular}{|l|l|*{5}{r}|*{5}{r}}
\hline
\multirow{2}{*}{Benchmark} & \multirow{2}{*}{Radius} & \multicolumn{5}{l|}{\texttt{RF}}                                                                                 	& \multicolumn{5}{l|}{\texttt{SFD}}                                      	\\ \cline{3-12}
                       	& & \multicolumn{1}{l|}{time} & \multicolumn{1}{l|}{$\sigma(\%)$} & \multicolumn{1}{l|}{$v_r(\%)$} &
                       	\multicolumn{1}{l|}{$v'_r(\%)$} &\multicolumn{1}{l|}{S/A/T} & \multicolumn{1}{l|}{time} & \multicolumn{1}{l|}{$\sigma(\%)$} & \multicolumn{1}{l|}{$v_r(\%)$} &
                       	\multicolumn{1}{l|}{$v'_r(\%)$} &\multicolumn{1}{l|}{S/A/T}  \\ \hline
\multirow{2}{*}{\texttt{city}}&Small & 1804.4 & 16.2\% & 54.1\% & 89.3\% & 0/9/9 & 0.2 & 0.3\% & 60.6\% & 100.0\% & \multicolumn{1}{l|}{9/9/9} \\
&Large & 1801.5 & 25.8\% & 21.3\% & 92.0\% & 0/9/9 & 3.4 & 1.1\% & 23.2\% & 100.0\% & \multicolumn{1}{l|}{9/9/9} \\ \cline{1-1}
\multirow{2}{*}{\texttt{Kgroup\_A}}&Small & 1622.6 & 21.5\% & 77.5\% & 91.1\% & 1/11/11 & 0.5 & 2.7\% & 85.1\% & 100.0\% & \multicolumn{1}{l|}{11/11/11} \\
&Large & 158.9 & 19.2\% & 19.3\% & 85.5\% & 8/11/11 & 0.4 & 6.7\% & 22.6\% & 100.0\% & \multicolumn{1}{l|}{11/11/11} \\ \cline{1-1}
\multirow{2}{*}{\texttt{Kgroup\_B}}&Small & 1800.9 & 59.1\% & 154.2\% & 185.0\% & 0/12/12 & 66.1 & 0.8\% & 83.3\% & 100.0\% & \multicolumn{1}{l|}{10/12/12} \\
&Large & 1800.6 & 75.5\% & 63.3\% & 312.1\% & 0/12/12 & 136.5 & 1.1\% & 20.2\% & 100.0\% & \multicolumn{1}{l|}{12/12/12} \\ \cline{1-1}
\multirow{2}{*}{\texttt{random\_A}}&Small & 15.9 & 8.1\% & 54.3\% & 86.0\% & 9/12/12 & 0.0 & 1.2\% & 63.2\% & 100.0\% & \multicolumn{1}{l|}{12/12/12} \\
&Large & 0.3 & 26.6\% & 19.8\% & 93.6\% & 12/12/12 & 0.0 & 2.4\% & 21.1\% & 100.0\% & \multicolumn{1}{l|}{12/12/12} \\ \cline{1-1}
\multirow{2}{*}{\texttt{random\_B}}&Small & 1304.3 & 38.5\% & 63.8\% & 91.8\% & 1/12/12 & 1.0 & 2.1\% & 69.5\% & 100.0\% & \multicolumn{1}{l|}{12/12/12} \\
&Large & 190.2 & 19.8\% & 11.2\% & 103.3\% & 9/12/12 & 1.7 & 8.7\% & 10.9\% & 100.0\% & \multicolumn{1}{l|}{12/12/12} \\ \cline{1-12}
\multirow{2}{*}{\texttt{all}}&Small & 604.9 & 23.7\% & 75.4\% & 104.6\% & 11/56/56 & 2.2 & 1.3\% & 72.1\% & 100.0\% & \multicolumn{1}{l|}{54/56/56} \\
&Large & 146.6 & 29.2\% & 22.8\% & 121.3\% & 29/56/56 & 3.9 & 3.2\% & 18.7\% & 100.0\% & \multicolumn{1}{l|}{56/56/56} \\
\cline{1-12}
\end{tabular}
\end{center}
\caption{Results for continuous and discrete models} \label{tab.cdm}
\end{table}

\Cref{tab.cdm} depicts some comparative results. A first observation is that \texttt{SFD} has fewer variables and constraints than  \texttt{RF}, as it models a simpler problem. In addition, our strengthening techniques explain that \texttt{SFD} can solve almost all instances in a very short time. Moreover, even the average relative primal bound of \texttt{SFD} is smaller than \texttt{RF}. However, we note that, with the exception of \texttt{Kgroup\_B} and \texttt{random\_B} of large radius, \texttt{RF} finds solutions with fewer facilities.  For \texttt{Kgroup\_B}, \texttt{RF} has a larger average dual gap than \texttt{SFD}.

In \Cref{fig.scatter2}, we  also show scatter plots of the relative dual gaps ($\sigma$) and the relative primal bounds ($v_r$) for those instances affected by both \texttt{SFD} and \texttt{RF}. These plots complement the averaged results of \Cref{tab.cdm} by giving information on all affected instances, and support the above analysis.

By allowing location at edges, the continuous model can reduce the number of installed facilities. However, it becomes more challenging to solve the problem. The results suggest that calling \texttt{SFD} and passing its solution as a warm-start to \texttt{RF} can make sense as a two-step optimization approach.

\section{Conclusions}\label{sec:conclu}
In this work, we use an integer programming approach to solve \problem,  propose various MILP formulations for this problem and test these formulations against an existing MILP formulation on several benchmarks from the literature.

The existing works mainly consider discretization methods and FDS. Discretization methods are indeed preprocessing procedures that restrict \problem\ to an equivalent set-covering problem with continuous demands and candidate facilities in FDS. But FDS is only computable when \problem\ satisfies some assumption, so it is not practical to employ FDS  for general \problem. On the other hand, to delimit the search space of MILP models, we use alternative preprocessing procedures. We explore the delimitation that relaxes the concept of FDS,  which also restricts the candidate space of facilities (from the full network to a still continuous sub-network).  We learn the following ideas to tackle similar problems (possibly with more complex constraints): i). Integer programming methods are more viable and flexible for general graphs, as a partial delimitation of the problems is useful for  strengthening the models (via separation of valid inequalities, tightening big-M constants, variable fixing), while the discretization methods require a full delimitation and a complete characterization of FDS. ii). Our delimitation is applicable as long as the facility location is continuous.

 Specifically, we devise and implement four models for \problem:   \texttt{F0},  \texttt{F} and  \texttt{SF}, which belong to the same family of models, and \texttt{RF}, which is the reduced one. These models mainly differ in preprocessing procedures applied. We find that  MILP solvers cannot read or build  the MILP model from \cite{Hamacher20} for any instances in our test bed, and this is due to the model having a large number of constraints and variables. So we mainly compare our four models. We find the MILP size is the main barrier to scalability. The delimitation of those parts of the network that can be covered from a specific location has been revealed as a very effective technique to reduce the model size. In addition, avoiding breaking long edges in the reduced model also results in better scalability. In conclusion,
\texttt{RF} is the best model: it can find good solutions with a small dual gap.

Meanwhile, the model  \texttt{SF}  is easily cast into   \texttt{SFD} for the discrete restriction of \problem. We find that allowing continuous facilities decreases the number of installed facilities but increases the solving time significantly.  We note that \texttt{SFD} finds an optimal solution for  the \dfproblem \, quickly, which is a primal solution for \problem. Therefore, \texttt{SFD} can be called as a fast MILP-based primal heuristic for \problem.

As for future studies, devising efficient heuristics to be integrated into MILP solvers can be useful to improve the primal performance of the proposed models. For instance, different relaxations of  \problem\ can be worth exploring, such as that where demand only happens at nodes (i.e. only nodes are to be covered). Every solution of \problem\, would be a solution of such relaxation, and hence the optimal value of the latter is a valid dual lower bound of the optimal value of  \problem.  If solving this combinatorial relaxation is efficient and provides a stronger dual lower bound than the LP relaxation of \problem, we can utilize this result and integrate the combinatorial dual bound into the MILP solver, which leads to a combinatorial branch-and-bound algorithm.

Moreover, we can investigate the potential use of FDS to further delimit the search space and integrate FDS in the preprocessing procedure.

\section*{Acknowledgments}
The authors would like to thank the anonymous reviewers for their careful reading of our manuscript and their many insightful comments and suggestions. The authors would also like to thank Joerg Kalcsics for providing access to some of the real-world network instances used for the experiments.

\section*{Funding source}
The research of Liding Xu is supported by French Ministry of Higher Education and Research (MESRI). The research of Mercedes Pelegr\'in was supported by the Chair “Integrated Urban Mobility”, backed by L’X -
École Polytechnique and La Fondation de l’École Polytechnique. The Partners of the Chair shall
not under any circumstances accept any liability for the content of this publication, for which the author shall be solely liable.

\bibliography{mybibfile}

\section*{Appendix}
\begin{landscape}
\scriptsize
\setlength{\tabcolsep}{1.9pt}
	\begin{longtable}{llrcrc|rrr|rrr|rrr|rrr|rrr|rrr|}
    	\caption{Detailed experimental results for comparing different models.
       	}
    	\label{table:root:detailed}\\
    	\tabledefline{SN}{subdivided network}
   	\tabledefline{DN}{degree-two-free network}
 	\tabledefline{time}{solving time}
	\tabledefline{gap}{relative dual gap $\sigma$}
	\tabledefline{primal}{primal bound $\underline{v}$}
	\tabledefline{-}{the result is not available (for the solver fails to load the MILP model)}
	\toprule
	& \multirow{2}{*}{Radius} & \multicolumn{2}{c}{\textbf{SN}} & \multicolumn{2}{c|}{\textbf{DN}} &   \multicolumn{3}{c|}{\texttt{EF}} & \multicolumn{3}{c|}{\texttt{F0}} & \multicolumn{3}{c|}{\texttt{F}} & \multicolumn{3}{c|}{\texttt{SF}} &
	\multicolumn{3}{c|}{\texttt{RF}} & \multicolumn{3}{c|}{\texttt{SFD}}  \\
        	instance & & nodes & edges & nodes   &  edges &  time & gap & primal & time & gap & primal& time & gap & primal & time & gap & primal & time & gap & primal & time & gap & primal \\
    	\toprule
    	\endfirsthead
    	\textbf{Table~\ref{table:root:detailed}} continued\\
    	\toprule
  	& \multirow{2}{*}{Radius} & \multicolumn{2}{c}{\textbf{SN}} & \multicolumn{2}{c|}{\textbf{DN}} &   \multicolumn{3}{c|}{\texttt{EF}} & \multicolumn{3}{c|}{\texttt{F0}} & \multicolumn{3}{c|}{\texttt{F}} & \multicolumn{3}{c|}{\texttt{SF}} &
	\multicolumn{3}{c|}{\texttt{RF}} & \multicolumn{3}{c|}{\texttt{SFD}}  \\
        	instance & & nodes & edges & nodes   &  edges &  time & gap & primal & time & gap & primal& time & gap & primal & time & gap & primal & time & gap & primal & time & gap & primal \\
    	\toprule
    	\endhead
\multirow{2}{*}{\texttt{city\_628}}&Small & 756 & 919 & 557 & 720&- & - & -&- & - & -&1800.1 & 42.0\% & 74.6\%&1801.7 & 57.8\% & 104.7\%&1803.5 & 18.3\% & 55.7\%&0.4 & 0.0\% & 60.6\%\\
&Large & 486 & 649 & 456 & 619&- & - & -&- & - & -&1800.7 & 28.0\% & 23.4\%&1800.6 & 31.5\% & 24.1\%&1801.8 & 24.6\% & 22.6\%&27.3 & 0.0\% & 23.4\%\\
\multirow{2}{*}{\texttt{city\_265}}&Small & 298 & 373 & 221 & 296&- & - & -&- & - & -&1800.4 & 25.9\% & 59.5\%&1802.7 & 26.6\% & 56.3\%&1801.0 & 16.9\% & 54.2\%&0.1 & 0.8\% & 62.1\%\\
&Large & 178 & 253 & 170 & 245&- & - & -&1801.4 & 28.3\% & 22.1\%&1801.3 & 28.6\% & 21.6\%&1801.1 & 32.0\% & 21.6\%&1802.4 & 26.5\% & 21.1\%&0.6 & 2.2\% & 24.2\%\\
\multirow{2}{*}{\texttt{city\_276}}&Small & 318 & 403 & 232 & 317&- & - & -&- & - & -&1800.1 & 41.9\% & 68.1\%&1800.8 & 65.4\% & 116.2\%&1805.9 & 21.8\% & 55.5\%&0.2 & 0.8\% & 61.8\%\\
&Large & 186 & 271 & 175 & 260&- & - & -&1802.3 & 33.3\% & 22.0\%&1800.6 & 39.3\% & 23.0\%&1800.8 & 37.4\% & 22.0\%&1800.6 & 34.4\% & 21.5\%&12.0 & 0.0\% & 22.5\%\\
\multirow{2}{*}{\texttt{city\_213}}&Small & 249 & 319 & 196 & 266&- & - & -&1811.4 & 28.8\% & 62.2\%&1800.7 & 30.4\% & 55.9\%&1800.3 & 30.6\% & 55.9\%&1807.9 & 24.3\% & 55.9\%&0.2 & 0.0\% & 61.5\%\\
&Large & 136 & 206 & 132 & 202&- & - & -&1800.7 & 33.5\% & 20.3\%&1800.5 & 44.6\% & 21.7\%&1801.0 & 42.1\% & 21.0\%&1800.8 & 41.4\% & 21.7\%&6.0 & 3.0\% & 23.1\%\\
\multirow{2}{*}{\texttt{city\_771}}&Small & 944 & 1123 & 755 & 934&- & - & -&- & - & -&1800.1 & 41.4\% & 72.3\%&1801.1 & 40.4\% & 71.3\%&1803.1 & 16.3\% & 54.6\%&0.3 & 0.0\% & 59.0\%\\
&Large & 591 & 770 & 565 & 744&- & - & -&- & - & -&1800.9 & 31.5\% & 22.3\%&1800.3 & 30.1\% & 21.8\%&1802.1 & 25.4\% & 20.9\%&6.3 & 0.0\% & 22.6\%\\
\multirow{2}{*}{\texttt{city\_138}}&Small & 162 & 201 & 135 & 174&- & - & -&1801.4 & 18.1\% & 55.6\%&1801.1 & 21.9\% & 54.5\%&1800.6 & 19.6\% & 54.5\%&1804.7 & 18.3\% & 52.5\%&0.0 & 0.0\% & 60.6\%\\
&Large & 88 & 127 & 88 & 127&- & - & -&1801.2 & 20.6\% & 20.2\%&1802.6 & 23.8\% & 20.2\%&1801.5 & 25.8\% & 20.2\%&1800.7 & 24.6\% & 20.2\%&0.3 & 4.3\% & 23.2\%\\
\multirow{2}{*}{\texttt{city\_479}}&Small & 584 & 689 & 428 & 533&- & - & -&- & - & -&1811.1 & 36.9\% & 69.3\%&1801.5 & 20.8\% & 56.1\%&1801.0 & 10.2\% & 52.1\%&0.2 & 0.5\% & 57.8\%\\
&Large & 379 & 484 & 352 & 457&- & - & -&- & - & -&1800.7 & 18.1\% & 21.7\%&1801.2 & 20.1\% & 22.7\%&1801.7 & 15.6\% & 21.7\%&4.5 & 1.1\% & 23.8\%\\
\multirow{2}{*}{\texttt{city\_132}}&Small & 152 & 179 & 112 & 139&- & - & -&1802.7 & 10.9\% & 55.7\%&1801.2 & 14.2\% & 53.8\%&1802.5 & 12.5\% & 52.8\%&1805.8 & 8.1\% & 51.9\%&0.0 & 1.5\% & 62.3\%\\
&Large & 91 & 118 & 87 & 114&- & - & -&1801.7 & 16.3\% & 20.8\%&1802.2 & 19.3\% & 20.8\%&1801.6 & 17.1\% & 20.8\%&1801.5 & 18.1\% & 20.8\%&0.0 & 4.0\% & 23.6\%\\
\multirow{2}{*}{\texttt{city\_268}}&Small & 307 & 380 & 229 & 302&- & - & -&- & - & -&1810.9 & 25.3\% & 56.4\%&1800.6 & 30.8\% & 60.5\%&1806.4 & 18.1\% & 54.9\%&0.1 & 0.0\% & 60.0\%\\
&Large & 188 & 261 & 178 & 251&- & - & -&1800.4 & 39.1\% & 24.1\%&1801.0 & 32.2\% & 20.5\%&1800.4 & 34.8\% & 21.5\%&1802.2 & 30.3\% & 21.5\%&1.5 & 0.0\% & 22.1\%\\
\multirow{2}{*}{\texttt{K100.5.con}}&Small & 86 & 175 & 64 & 153&- & - & -&1811.1 & 34.5\% & 80.6\%&1800.4 & 40.6\% & 77.4\%&1802.6 & 46.9\% & 80.6\%&1802.9 & 33.4\% & 77.4\%&0.4 & 4.2\% & 77.4\%\\
&Large & 39 & 128 & 37 & 126&- & - & -&1800.4 & 25.0\% & 25.8\%&1555.1 & 14.3\% & 22.6\%&1395.3 & 14.3\% & 22.6\%&1379.2 & 14.3\% & 22.6\%&0.1 & 0.0\% & 25.8\%\\
\multirow{2}{*}{\texttt{K100.2}}&Small & 61 & 120 & 42 & 101&- & - & -&1800.9 & 23.0\% & 79.2\%&1800.6 & 29.8\% & 79.2\%&1802.2 & 24.4\% & 75.0\%&1800.7 & 18.8\% & 75.0\%&0.1 & 4.8\% & 87.5\%\\
&Large & 32 & 91 & 32 & 91&- & - & -&42.5 & 25.0\% & 16.7\%&12.8 & 25.0\% & 16.7\%&74.4 & 25.0\% & 16.7\%&90.9 & 25.0\% & 16.7\%&0.1 & 20.0\% & 20.8\%\\
\multirow{2}{*}{\texttt{K100.3.con}}&Small & 94 & 191 & 66 & 163&- & - & -&1800.7 & 41.2\% & 103.8\%&1805.3 & 50.0\% & 100.0\%&1801.2 & 46.2\% & 96.2\%&1800.5 & 26.6\% & 92.3\%&3.4 & 4.2\% & 92.3\%\\
&Large & 40 & 137 & 29 & 126&- & - & -&436.3 & 0.0\% & 26.9\%&425.4 & 14.3\% & 26.9\%&409.3 & 14.3\% & 26.9\%&22.2 & 14.3\% & 26.9\%&0.2 & 12.5\% & 30.8\%\\
\multirow{2}{*}{\texttt{K100.10}}&Small & 67 & 118 & 51 & 102&- & - & -&1800.8 & 12.0\% & 70.4\%&1800.4 & 23.3\% & 74.1\%&1801.4 & 27.4\% & 74.1\%&1800.9 & 19.8\% & 70.4\%&0.1 & 4.2\% & 88.9\%\\
&Large & 35 & 86 & 35 & 86&- & - & -&91.5 & 20.0\% & 18.5\%&163.5 & 20.0\% & 18.5\%&328.2 & 20.0\% & 18.5\%&164.5 & 20.0\% & 18.5\%&0.1 & 16.7\% & 22.2\%\\
\multirow{2}{*}{\texttt{K100.con}}&Small & 145 & 291 & 121 & 267&- & - & -&1800.2 & 49.8\% & 86.7\%&1806.8 & 55.7\% & 86.7\%&1800.2 & 57.8\% & 88.9\%&1800.1 & 46.3\% & 77.8\%&0.6 & 0.0\% & 80.0\%\\
&Large & 55 & 201 & 55 & 201&- & - & -&1800.9 & 36.7\% & 20.0\%&1802.9 & 36.6\% & 20.0\%&1802.3 & 36.5\% & 20.0\%&1802.2 & 44.3\% & 22.2\%&1.8 & 0.0\% & 22.2\%\\
\multirow{2}{*}{\texttt{K100.4.con}}&Small & 85 & 169 & 60 & 144&- & - & -&1806.5 & 34.8\% & 93.1\%&1801.7 & 37.0\% & 86.2\%&1800.7 & 23.1\% & 75.9\%&1804.5 & 23.4\% & 75.9\%&0.2 & 4.2\% & 82.8\%\\
&Large & 40 & 124 & 40 & 124&- & - & -&187.2 & 0.0\% & 13.8\%&254.1 & 0.0\% & 13.8\%&1723.2 & 0.0\% & 13.8\%&92.7 & 0.0\% & 13.8\%&0.2 & 0.0\% & 17.2\%\\
\multirow{2}{*}{\texttt{K100.9}}&Small & 56 & 104 & 35 & 83&- & - & -&1803.0 & 9.3\% & 81.8\%&1805.8 & 18.7\% & 72.7\%&1800.4 & 12.2\% & 68.2\%&565.2 & 6.7\% & 68.2\%&0.1 & 5.6\% & 81.8\%\\
&Large & 30 & 78 & 30 & 78&- & - & -&1.3 & 33.3\% & 13.6\%&3.6 & 33.3\% & 13.6\%&2.4 & 33.3\% & 13.6\%&2.6 & 33.3\% & 13.6\%&0.1 & 25.0\% & 18.2\%\\
\multirow{2}{*}{\texttt{K100.1}}&Small & 120 & 263 & 70 & 213&- & - & -&1800.2 & 37.4\% & 109.5\%&1806.6 & 45.3\% & 100.0\%&1800.1 & 52.6\% & 97.6\%&1810.0 & 16.7\% & 85.7\%&1.9 & 2.4\% & 97.6\%\\
&Large & 65 & 208 & 61 & 204&- & - & -&1801.7 & 34.9\% & 23.8\%&1800.3 & 43.9\% & 26.2\%&1800.2 & 45.5\% & 26.2\%&1800.4 & 44.6\% & 26.2\%&1.0 & 9.1\% & 26.2\%\\
\multirow{2}{*}{\texttt{K100.7}}&Small & 74 & 142 & 48 & 116&- & - & -&1800.9 & 24.0\% & 92.0\%&1803.7 & 31.3\% & 92.0\%&1800.4 & 39.2\% & 92.0\%&1801.7 & 22.5\% & 92.0\%&0.1 & 4.0\% & 100.0\%\\
&Large & 37 & 105 & 37 & 105&- & - & -&7.9 & 25.0\% & 16.0\%&214.0 & 25.0\% & 16.0\%&257.0 & 25.0\% & 16.0\%&141.9 & 25.0\% & 16.0\%&0.1 & 20.0\% & 20.0\%\\
\multirow{2}{*}{\texttt{K100.6}}&Small & 50 & 92 & 38 & 80&- & - & -&1803.6 & 11.1\% & 68.2\%&1800.9 & 20.1\% & 68.2\%&1804.9 & 21.9\% & 68.2\%&1807.3 & 16.3\% & 68.2\%&0.1 & 0.0\% & 72.7\%\\
&Large & 28 & 70 & 28 & 70&- & - & -&4.6 & 25.0\% & 18.2\%&64.0 & 25.0\% & 18.2\%&14.1 & 25.0\% & 18.2\%&11.7 & 25.0\% & 18.2\%&0.0 & 20.0\% & 22.7\%\\
\multirow{2}{*}{\texttt{K100.8.con}}&Small & 107 & 208 & 78 & 179&- & - & -&1800.6 & 36.2\% & 79.1\%&1800.8 & 33.4\% & 74.4\%&1800.5 & 32.9\% & 76.7\%&1800.7 & 29.0\% & 74.4\%&0.9 & 2.9\% & 79.1\%\\
&Large & 57 & 158 & 53 & 154&- & - & -&1803.0 & 30.0\% & 23.3\%&1800.4 & 38.3\% & 23.3\%&1800.5 & 39.6\% & 23.3\%&1800.9 & 36.4\% & 23.3\%&2.9 & 9.1\% & 25.6\%\\
\multirow{2}{*}{\texttt{K400.9.con}}&Small & 588 & 1239 & 363 & 1014&- & - & -&- & - & -&1800.1 & 87.4\% & 278.2\%&1803.5 & 84.1\% & 216.6\%&1801.5 & 36.2\% & 93.8\%&23.8 & 0.5\% & 88.2\%\\
&Large & 312 & 963 & 284 & 935&- & - & -&- & - & -&1800.3 & 84.9\% & 53.6\%&1800.8 & 93.1\% & 116.6\%&1800.4 & 70.6\% & 42.7\%&81.6 & 0.0\% & 22.3\%\\
\multirow{2}{*}{\texttt{K400.3.con}}&Small & 595 & 1191 & 450 & 1046&- & - & -&- & - & -&1801.4 & 89.8\% & 301.9\%&1801.0 & 85.7\% & 214.8\%&1800.1 & 56.4\% & 101.4\%&96.1 & 0.6\% & 81.4\%\\
&Large & 273 & 869 & 249 & 845&- & - & -&- & - & -&1800.5 & 84.1\% & 59.0\%&1800.3 & 81.7\% & 49.5\%&1800.5 & 62.6\% & 33.8\%&96.6 & 0.0\% & 21.9\%\\
\multirow{2}{*}{\texttt{K200.con}}&Small & 184 & 374 & 121 & 311&- & - & -&1804.4 & 39.3\% & 86.4\%&1800.4 & 39.0\% & 80.2\%&1805.2 & 45.5\% & 74.1\%&1800.7 & 27.9\% & 72.8\%&0.6 & 1.5\% & 80.2\%\\
&Large & 109 & 299 & 109 & 299&- & - & -&1801.3 & 42.7\% & 17.3\%&1800.6 & 48.4\% & 18.5\%&1800.3 & 57.2\% & 22.2\%&1801.3 & 47.7\% & 18.5\%&10.4 & 6.7\% & 18.5\%\\
\multirow{2}{*}{\texttt{K400}}&Small & 715 & 1398 & 568 & 1251&- & - & -&- & - & -&1801.2 & 90.1\% & 310.0\%&1800.2 & 86.8\% & 237.7\%&1800.1 & 86.5\% & 315.6\%&371.9 & 0.5\% & 79.2\%\\
&Large & 296 & 979 & 284 & 967&- & - & -&- & - & -&1800.4 & 93.4\% & 125.5\%&1800.3 & 92.8\% & 112.6\%&1800.9 & 78.8\% & 47.6\%&211.2 & 0.0\% & 18.6\%\\
\multirow{2}{*}{\texttt{K400.4}}&Small & 516 & 1103 & 358 & 945&- & - & -&- & - & -&1800.1 & 87.4\% & 255.3\%&1802.3 & 77.6\% & 138.6\%&1800.7 & 51.0\% & 105.1\%&17.3 & 0.6\% & 87.3\%\\
&Large & 270 & 857 & 256 & 843&- & - & -&- & - & -&1800.6 & 94.5\% & 142.1\%&1801.1 & 93.3\% & 115.2\%&1800.2 & 67.0\% & 29.9\%&74.2 & 0.0\% & 19.8\%\\
\multirow{2}{*}{\texttt{K400.10.con}}&Small & 671 & 1373 & 486 & 1188&- & - & -&- & - & -&1801.7 & 89.2\% & 294.1\%&1800.7 & 86.3\% & 230.3\%&1800.1 & 73.9\% & 194.1\%&27.4 & 0.0\% & 83.3\%\\
&Large & 300 & 1002 & 291 & 993&- & - & -&- & - & -&1800.1 & 89.0\% & 63.8\%&1800.7 & 94.2\% & 119.5\%&1800.2 & 94.0\% & 130.8\%&116.4 & 2.4\% & 19.0\%\\
\multirow{2}{*}{\texttt{K400.2}}&Small & 709 & 1429 & 522 & 1242&- & - & -&- & - & -&1800.6 & 89.8\% & 321.1\%&1800.2 & 87.2\% & 247.8\%&1801.3 & 52.4\% & 107.9\%&32.4 & 0.5\% & 86.4\%\\
&Large & 311 & 1031 & 305 & 1025&- & - & -&- & - & -&1800.3 & 94.2\% & 139.9\%&1801.6 & 84.4\% & 49.6\%&1800.3 & 94.1\% & 143.4\%&221.0 & 2.4\% & 18.0\%\\
\multirow{2}{*}{\texttt{K400.6.con}}&Small & 789 & 1583 & 612 & 1406&- & - & -&- & - & -&1800.0 & 91.7\% & 347.3\%&1803.0 & 88.8\% & 255.2\%&1800.6 & 86.7\% & 346.9\%&1801.0 & 2.2\% & 81.7\%\\
&Large & 321 & 1115 & 301 & 1095&- & - & -&- & - & -&1800.5 & 93.7\% & 129.9\%&1800.9 & 93.2\% & 118.3\%&1800.7 & 92.1\% & 132.0\%&399.0 & 2.0\% & 20.3\%\\
\multirow{2}{*}{\texttt{K400.1.con}}&Small & 587 & 1224 & 379 & 1016&- & - & -&- & - & -&1800.5 & 69.3\% & 155.8\%&1800.7 & 76.0\% & 203.7\%&1800.5 & 65.3\% & 173.3\%&15.6 & 0.5\% & 89.9\%\\
&Large & 314 & 951 & 259 & 896&- & - & -&- & - & -&1800.2 & 84.6\% & 78.8\%&1800.3 & 87.1\% & 115.7\%&1801.4 & 51.8\% & 40.1\%&82.8 & 0.0\% & 26.3\%\\
\multirow{2}{*}{\texttt{K400.8.con}}&Small & 749 & 1501 & 610 & 1362&- & - & -&- & - & -&1800.0 & 91.6\% & 337.9\%&1800.0 & 88.8\% & 246.0\%&1801.2 & 88.7\% & 347.2\%&1800.3 & 4.6\% & 77.0\%\\
&Large & 294 & 1046 & 284 & 1036&- & - & -&- & - & -&1800.6 & 93.8\% & 123.0\%&1800.9 & 93.4\% & 116.6\%&1800.1 & 92.8\% & 122.1\%&1183.8 & 2.2\% & 19.6\%\\
\multirow{2}{*}{\texttt{K400.7.con}}&Small & 633 & 1275 & 462 & 1104&- & - & -&- & - & -&1801.2 & 88.5\% & 280.0\%&1800.2 & 85.7\% & 217.3\%&1802.2 & 72.3\% & 178.7\%&140.4 & 0.5\% & 83.6\%\\
&Large & 302 & 944 & 290 & 932&- & - & -&- & - & -&1800.3 & 80.1\% & 44.4\%&1800.5 & 86.0\% & 62.2\%&1800.6 & 84.0\% & 62.7\%&264.0 & 2.3\% & 19.6\%\\
\multirow{2}{*}{\texttt{K400.5.con}}&Small & 600 & 1179 & 431 & 1010&- & - & -&- & - & -&1800.3 & 77.9\% & 144.1\%&1800.2 & 76.9\% & 137.7\%&1801.4 & 50.6\% & 99.5\%&43.6 & 0.0\% & 82.7\%\\
&Large & 295 & 874 & 285 & 864&- & - & -&- & - & -&1800.2 & 94.2\% & 136.4\%&1800.3 & 82.6\% & 45.0\%&1800.4 & 92.9\% & 135.9\%&112.4 & 2.3\% & 20.0\%\\
\multirow{2}{*}{\texttt{r\_15\_0.3\_25}}&Small & 26 & 36 & 21 & 31&- & - & -&20.2 & 12.5\% & 53.3\%&56.8 & 12.5\% & 53.3\%&23.9 & 12.5\% & 53.3\%&13.6 & 12.5\% & 53.3\%&0.0 & 0.0\% & 60.0\%\\
&Large & 15 & 25 & 15 & 25&- & - & -&0.1 & 33.3\% & 20.0\%&0.4 & 33.3\% & 20.0\%&0.2 & 33.3\% & 20.0\%&0.4 & 33.3\% & 20.0\%&0.0 & 33.3\% & 20.0\%\\
\multirow{2}{*}{\texttt{r\_10\_0.3\_12}}&Small & 15 & 17 & 13 & 15&- & - & -&0.3 & 16.7\% & 60.0\%&0.2 & 16.7\% & 60.0\%&0.1 & 16.7\% & 60.0\%&0.2 & 16.7\% & 60.0\%&0.0 & 14.3\% & 70.0\%\\
&Large & 10 & 12 & 10 & 12&- & - & -&0.0 & 33.3\% & 30.0\%&0.0 & 33.3\% & 30.0\%&0.0 & 33.3\% & 30.0\%&0.0 & 33.3\% & 30.0\%&0.0 & 0.0\% & 30.0\%\\
\multirow{2}{*}{\texttt{r\_20\_0.2\_34}}&Small & 39 & 53 & 33 & 47&- & - & -&336.8 & 9.1\% & 55.0\%&376.5 & 9.1\% & 55.0\%&159.1 & 9.1\% & 55.0\%&178.8 & 9.1\% & 55.0\%&0.0 & 0.0\% & 60.0\%\\
&Large & 20 & 34 & 20 & 34&- & - & -&0.2 & 33.3\% & 15.0\%&0.2 & 33.3\% & 15.0\%&0.2 & 33.3\% & 15.0\%&0.3 & 33.3\% & 15.0\%&0.0 & 0.0\% & 15.0\%\\
\multirow{2}{*}{\texttt{r\_15\_0.4\_45}}&Small & 40 & 70 & 40 & 70&- & - & -&1800.7 & 20.7\% & 53.3\%&1808.2 & 23.6\% & 53.3\%&1801.2 & 37.4\% & 60.0\%&1800.9 & 28.4\% & 53.3\%&0.1 & 10.0\% & 66.7\%\\
&Large & 15 & 45 & 15 & 45&- & - & -&0.2 & 50.0\% & 13.3\%&0.6 & 50.0\% & 13.3\%&0.6 & 50.0\% & 13.3\%&0.7 & 50.0\% & 13.3\%&0.1 & 50.0\% & 13.3\%\\
\multirow{2}{*}{\texttt{r\_10\_0.1\_12}}&Small & 16 & 18 & 14 & 16&- & - & -&0.1 & 20.0\% & 50.0\%&0.1 & 20.0\% & 50.0\%&0.0 & 20.0\% & 50.0\%&0.0 & 20.0\% & 50.0\%&0.0 & 0.0\% & 60.0\%\\
&Large & 10 & 12 & 10 & 12&- & - & -&0.0 & 33.3\% & 30.0\%&0.0 & 33.3\% & 30.0\%&0.0 & 0.0\% & 20.0\%&0.0 & 33.3\% & 30.0\%&0.0 & 0.0\% & 20.0\%\\
\multirow{2}{*}{\texttt{r\_20\_0.3\_49}}&Small & 46 & 75 & 46 & 75&- & - & -&1801.4 & 18.2\% & 55.0\%&1801.6 & 22.7\% & 55.0\%&1800.6 & 21.8\% & 55.0\%&1800.7 & 22.6\% & 55.0\%&0.0 & 7.1\% & 70.0\%\\
&Large & 20 & 49 & 20 & 49&- & - & -&0.3 & 50.0\% & 10.0\%&0.3 & 50.0\% & 10.0\%&1.0 & 50.0\% & 10.0\%&1.1 & 50.0\% & 10.0\%&0.1 & 33.3\% & 15.0\%\\
\multirow{2}{*}{\texttt{r\_20\_0.1\_23}}&Small & 29 & 32 & 23 & 26&- & - & -&2.1 & 9.1\% & 55.0\%&2.8 & 10.0\% & 50.0\%&1.0 & 10.0\% & 50.0\%&2.4 & 0.0\% & 50.0\%&0.0 & 0.0\% & 60.0\%\\
&Large & 20 & 23 & 20 & 23&- & - & -&0.0 & 20.0\% & 25.0\%&0.2 & 20.0\% & 25.0\%&0.1 & 0.0\% & 20.0\%&0.1 & 0.0\% & 20.0\%&0.0 & 0.0\% & 25.0\%\\
\multirow{2}{*}{\texttt{r\_15\_0.1\_14}}&Small & 19 & 18 & 15 & 14&- & - & -&0.0 & 12.5\% & 53.3\%&0.0 & 12.5\% & 53.3\%&0.0 & 12.5\% & 53.3\%&0.0 & 0.0\% & 53.3\%&0.0 & 0.0\% & 53.3\%\\
&Large & 15 & 14 & 15 & 14&- & - & -&0.0 & 25.0\% & 26.7\%&0.0 & 25.0\% & 26.7\%&0.0 & 25.0\% & 26.7\%&0.0 & 25.0\% & 26.7\%&0.0 & 0.0\% & 26.7\%\\
\multirow{2}{*}{\texttt{r\_20\_0.4\_69}}&Small & 52 & 101 & 52 & 101&- & - & -&1801.3 & 37.0\% & 60.0\%&1801.0 & 36.4\% & 60.0\%&1803.3 & 36.1\% & 60.0\%&1801.1 & 37.6\% & 60.0\%&0.1 & 6.7\% & 75.0\%\\
&Large & 20 & 69 & 20 & 69&- & - & -&1.2 & 50.0\% & 10.0\%&3.9 & 50.0\% & 10.0\%&1.0 & 50.0\% & 10.0\%&2.5 & 50.0\% & 10.0\%&0.2 & 33.3\% & 15.0\%\\
\multirow{2}{*}{\texttt{r\_15\_0.2\_22}}&Small & 25 & 32 & 21 & 28&- & - & -&3.4 & 12.5\% & 53.3\%&4.6 & 12.5\% & 53.3\%&7.0 & 12.5\% & 53.3\%&3.2 & 12.5\% & 53.3\%&0.0 & 0.0\% & 66.7\%\\
&Large & 15 & 22 & 15 & 22&- & - & -&0.0 & 33.3\% & 20.0\%&0.1 & 33.3\% & 20.0\%&0.1 & 33.3\% & 20.0\%&0.1 & 33.3\% & 20.0\%&0.0 & 0.0\% & 26.7\%\\
\multirow{2}{*}{\texttt{r\_10\_0.4\_13}}&Small & 14 & 17 & 14 & 17&- & - & -&0.2 & 16.7\% & 60.0\%&1.2 & 16.7\% & 60.0\%&0.4 & 16.7\% & 60.0\%&1.0 & 16.7\% & 60.0\%&0.0 & 0.0\% & 60.0\%\\
&Large & 10 & 13 & 10 & 13&- & - & -&0.0 & 33.3\% & 30.0\%&0.0 & 0.0\% & 20.0\%&0.0 & 0.0\% & 20.0\%&0.0 & 33.3\% & 30.0\%&0.0 & 0.0\% & 20.0\%\\
\multirow{2}{*}{\texttt{r\_10\_0.2\_9}}&Small & 12 & 11 & 5 & 4&- & - & -&0.0 & 20.0\% & 50.0\%&0.0 & 20.0\% & 50.0\%&0.0 & 20.0\% & 50.0\%&0.0 & 0.0\% & 50.0\%&0.0 & 0.0\% & 60.0\%\\
&Large & 10 & 9 & 10 & 9&- & - & -&0.0 & 0.0\% & 20.0\%&0.0 & 33.3\% & 30.0\%&0.0 & 0.0\% & 20.0\%&0.0 & 33.3\% & 30.0\%&0.0 & 0.0\% & 40.0\%\\
\multirow{2}{*}{\texttt{r\_40\_0.2\_148}}&Small & 112 & 220 & 112 & 220&- & - & -&1800.3 & 52.6\% & 65.0\%&1800.6 & 58.7\% & 67.5\%&1800.2 & 64.8\% & 72.5\%&1802.5 & 59.2\% & 67.5\%&0.5 & 3.4\% & 72.5\%\\
&Large & 40 & 148 & 40 & 148&- & - & -&1800.4 & 50.0\% & 10.0\%&1801.8 & 53.9\% & 10.0\%&1800.5 & 62.3\% & 10.0\%&1800.5 & 50.0\% & 10.0\%&7.2 & 25.0\% & 10.0\%\\
\multirow{2}{*}{\texttt{r\_40\_0.4\_297}}&Small & 191 & 448 & 191 & 448&- & - & -&1800.5 & 75.5\% & 80.0\%&1805.8 & 79.8\% & 92.5\%&1800.1 & 88.8\% & 150.0\%&1800.2 & 75.9\% & 77.5\%&22.2 & 3.4\% & 72.5\%\\
&Large & 40 & 297 & 40 & 297&- & - & -&1212.6 & 50.0\% & 5.0\%&1801.3 & 66.7\% & 7.5\%&1197.1 & 50.0\% & 5.0\%&1801.3 & 75.0\% & 10.0\%&7.7 & 0.0\% & 5.0\%\\
\multirow{2}{*}{\texttt{r\_30\_0.2\_75}}&Small & 67 & 112 & 67 & 112&- & - & -&1800.7 & 33.1\% & 60.0\%&1801.0 & 34.1\% & 60.0\%&1800.4 & 35.8\% & 60.0\%&1801.6 & 29.4\% & 56.7\%&0.2 & 0.0\% & 70.0\%\\
&Large & 30 & 75 & 30 & 75&- & - & -&844.2 & 25.0\% & 13.3\%&435.7 & 25.0\% & 13.3\%&854.9 & 25.0\% & 13.3\%&603.6 & 25.0\% & 13.3\%&1.6 & 20.0\% & 16.7\%\\
\multirow{2}{*}{\texttt{r\_25\_0.1\_36}}&Small & 40 & 51 & 36 & 47&- & - & -&41.3 & 7.7\% & 52.0\%&354.0 & 7.7\% & 52.0\%&201.9 & 7.7\% & 52.0\%&36.4 & 7.7\% & 52.0\%&0.0 & 0.0\% & 60.0\%\\
&Large & 25 & 36 & 25 & 36&- & - & -&6.2 & 0.0\% & 16.0\%&1.5 & 0.0\% & 16.0\%&0.7 & 0.0\% & 16.0\%&1.6 & 0.0\% & 16.0\%&0.0 & 0.0\% & 20.0\%\\
\multirow{2}{*}{\texttt{r\_25\_0.4\_112}}&Small & 83 & 170 & 83 & 170&- & - & -&1800.1 & 64.7\% & 72.0\%&1806.2 & 61.8\% & 64.0\%&1801.0 & 64.2\% & 68.0\%&1800.2 & 61.5\% & 64.0\%&0.3 & 5.9\% & 68.0\%\\
&Large & 25 & 112 & 25 & 112&- & - & -&0.3 & 50.0\% & 8.0\%&42.1 & 50.0\% & 8.0\%&53.2 & 50.0\% & 8.0\%&69.4 & 50.0\% & 8.0\%&0.5 & 50.0\% & 8.0\%\\
\multirow{2}{*}{\texttt{r\_25\_0.3\_98}}&Small & 74 & 147 & 74 & 147&- & - & -&1800.7 & 49.7\% & 64.0\%&1811.2 & 50.0\% & 64.0\%&1800.4 & 56.6\% & 68.0\%&1804.9 & 51.9\% & 64.0\%&0.2 & 5.3\% & 76.0\%\\
&Large & 25 & 98 & 25 & 98&- & - & -&86.2 & 50.0\% & 8.0\%&115.7 & 0.0\% & 8.0\%&40.6 & 0.0\% & 8.0\%&32.8 & 0.0\% & 8.0\%&1.0 & 33.3\% & 12.0\%\\
\multirow{2}{*}{\texttt{r\_40\_0.1\_84}}&Small & 78 & 122 & 76 & 120&- & - & -&1801.4 & 22.4\% & 52.5\%&1801.5 & 22.6\% & 52.5\%&1800.7 & 21.0\% & 52.5\%&1801.5 & 25.2\% & 55.0\%&0.1 & 0.0\% & 67.5\%\\
&Large & 40 & 84 & 40 & 84&- & - & -&104.7 & 16.7\% & 15.0\%&120.7 & 16.7\% & 15.0\%&61.8 & 16.7\% & 15.0\%&19.1 & 16.7\% & 15.0\%&0.4 & 14.3\% & 17.5\%\\
\multirow{2}{*}{\texttt{r\_25\_0.2\_58}}&Small & 53 & 86 & 53 & 86&- & - & -&1801.2 & 19.0\% & 52.0\%&1800.8 & 23.1\% & 52.0\%&1801.9 & 22.8\% & 52.0\%&1800.8 & 27.1\% & 56.0\%&0.1 & 5.9\% & 68.0\%\\
&Large & 25 & 58 & 25 & 58&- & - & -&163.3 & 25.0\% & 16.0\%&178.1 & 25.0\% & 16.0\%&196.7 & 25.0\% & 16.0\%&262.7 & 25.0\% & 16.0\%&0.3 & 25.0\% & 16.0\%\\
\multirow{2}{*}{\texttt{r\_40\_0.3\_219}}&Small & 149 & 328 & 149 & 328&- & - & -&1809.4 & 68.8\% & 77.5\%&1800.7 & 70.1\% & 75.0\%&1802.9 & 72.5\% & 72.5\%&1806.4 & 71.0\% & 77.5\%&3.1 & 3.7\% & 67.5\%\\
&Large & 40 & 219 & 40 & 219&- & - & -&823.1 & 50.0\% & 5.0\%&1642.4 & 50.0\% & 5.0\%&724.6 & 50.0\% & 5.0\%&1802.7 & 75.0\% & 10.0\%&8.2 & 50.0\% & 5.0\%\\
\multirow{2}{*}{\texttt{r\_30\_0.3\_131}}&Small & 92 & 193 & 92 & 193&- & - & -&1800.8 & 61.1\% & 66.7\%&1800.1 & 67.4\% & 73.3\%&1800.3 & 65.6\% & 66.7\%&1800.5 & 67.3\% & 73.3\%&0.7 & 4.5\% & 73.3\%\\
&Large & 30 & 131 & 30 & 131&- & - & -&1588.5 & 33.3\% & 10.0\%&1800.7 & 54.5\% & 10.0\%&894.2 & 33.3\% & 10.0\%&502.2 & 33.3\% & 10.0\%&3.2 & 33.3\% & 10.0\%\\
\multirow{2}{*}{\texttt{r\_30\_0.4\_188}}&Small & 117 & 275 & 117 & 275&- & - & -&1800.4 & 70.1\% & 73.3\%&1800.2 & 73.5\% & 86.7\%&1801.7 & 73.2\% & 80.0\%&1800.7 & 71.6\% & 76.7\%&5.4 & 4.3\% & 76.7\%\\
&Large & 30 & 188 & 30 & 188&- & - & -&99.8 & 50.0\% & 6.7\%&225.3 & 50.0\% & 6.7\%&328.3 & 50.0\% & 6.7\%&450.9 & 50.0\% & 6.7\%&2.2 & 0.0\% & 6.7\%\\
\multirow{2}{*}{\texttt{r\_30\_0.1\_54}}&Small & 52 & 76 & 46 & 70&- & - & -&1801.6 & 11.7\% & 53.3\%&1801.9 & 15.2\% & 53.3\%&1802.4 & 15.4\% & 53.3\%&1802.4 & 12.5\% & 53.3\%&0.0 & 0.0\% & 63.3\%\\
&Large & 30 & 54 & 30 & 54&- & - & -&57.8 & 20.0\% & 16.7\%&53.2 & 20.0\% & 16.7\%&48.4 & 20.0\% & 16.7\%&89.4 & 20.0\% & 16.7\%&0.2 & 0.0\% & 16.7\%\\
\bottomrule
\end{longtable}
\end{landscape}
\end{document}